\documentclass[12pt,twoside,a4paper,notitlepage]{article}

\usepackage{fancyhdr}

\usepackage{amsmath}

\usepackage{amssymb}
\usepackage{amsfonts}

\usepackage[a4paper,margin=2cm,vmargin={2cm,2cm},includeheadfoot]{geometry}
\usepackage{amsthm}
\usepackage{enumerate}
\usepackage[scaled=0.92]{helvet}

\usepackage{graphicx}
\usepackage[labelfont=bf,margin=1cm,justification=justified,labelsep=period]{caption}
\usepackage{subfigure}

\theoremstyle{plain}
\newtheorem{theorem}{Theorem}[section]
\newtheorem{proposition}[theorem]{Proposition}

\newtheorem{lemma}[theorem]{Lemma}

{\theoremstyle{definition}
\newtheorem{remark}[theorem]{Remark}}

\newcounter{appa}
\newtheorem{propositiona}[appa]{Proposition A.\!}

\newcounter{appb}

\newtheorem{theoremb}[appb]{Theorem B.\!}
\newtheorem{lemmab}[appb]{Lemma B.\!}

{\theoremstyle{definition}
\newtheorem{remarkb}[appb]{Remark B.\!}}

\newcommand{\R}{{\mathbb {R}}}  
\newcommand{\C}{{\mathbb{C}}}   %
\newcommand{\Z}{{\mathbb{Z}}}   %
\newcommand{\Q}{{\mathbb {Q}}}
\newcommand{\E}{{\mathbb{E}}}   
\newcommand{\Prob}{\mathbb{P}}  

\newcommand{\bs}{\boldsymbol}   
\newcommand{\wt}{\widetilde}    
\newcommand{\1}{{\bf 1}}        

\newcommand{\itn}{\it \bfseries  }

\newcommand{\inte}{\mathrel{\mathop{\kern0pt D}\limits^{\raisebox{-5pt}{\text{\scriptsize$\circ$}}}}}

\makeatletter
\renewcommand\@makefntext
{\noindent}
\makeatother

\begin{document}

\begin{center}

{\LARGE \bf A new look at the  Heston characteristic function} \\ [.5cm]

{\large \bf Sebastian del Ba\~{n}o Rollin,  Albert Ferreiro-Castilla,
  Frederic Utzet}

\end{center}
\thispagestyle{empty}


\footnotetext{S. del Ba{\~n}o Rollin, Centre de Recerca Matem{\`a}tica, Apartat 50, 08193 Bellaterra (Barcelona) Spain.
 e-mail:
{\tt sebastiandb@crm.cat}}

\footnotetext{A. Ferreiro-Castilla, Departament de Matem{\`a}tiques, Edifici C,
 Universitat Aut{\`o}noma de Barcelona,   08193 Bellaterra (Barcelona) Spain.  e-mail: {\tt aferreiro@mat.uab.cat}.
Supported by a  PhD grant of the   Centre de Recerca Matem{\`a}tica, and by  grant BFM2006-06427
 Ministerio de Educaci{\'o}n y Ciencia and FEDER}

 \footnotetext{F. Utzet (corresponding author), Departament de Matem{\`a}tiques, Edifici C,
 Universitat Aut{\`o}noma de Barcelona,   08193 Bellaterra (Barcelona) Spain.  e-mail: {\tt utzet@mat.uab.cat}.
 Supported by grant BFM2006-06427 Ministerio de Educaci{\'o}n y Ciencia and FEDER}

\noindent{\bf Abstract}
A new expression for the characteristic function of log-spot in Heston model
is presented.
This expression more clearly exhibits its properties as an analytic
characteristic function and allows
us to compute the exact domain of the moment generating function. This
result is then applied to the
volatility smile at extreme strikes and to the control of the moments of
spot. We also give a factorization of the moment generating function as
product of Bessel type factors, and an approximating sequence to the law of
log-spot is deduced.

\bigskip

\noindent{\bf Keywords} Heston volatility model, Characteristic function, Extreme strikes, Bessel random variables.

\bigskip

\noindent{\bf Mathematics Subject Classification (2000)} 91B28, 60H10, 60E10

\bigskip

\noindent{\bf JEL Classification} G13, C65

\section{Introduction}

The first surprising fact about the Heston stochastic  volatility model (Heston \cite{Hes93}) is that
the characteristic function of log-spot is
computable and has a nice expression in terms of elementary functions; its deduction
 had
enormous merit. The second thing, and still more fascinating, is that such characteristic function is analytic,
 that means (see Lukacs \cite[chapter 7]{Luk70}  for an equivalent   definition and the main properties of
  analytic characteristic functions)
there is a  function $\Psi(z)$ of the complex variable $z$, analytic in a neighborhood of 0, such that
$$\E[e^{i u X_t}]=\Psi(iu)$$
for $u$ (real) in a neighborhood of 0, where $X_t$ is  log-spot at time $t$.

The fact that a characteristic function is analytic has important consequences, and we here   exploit some of them.
However, the standard form of writing  the Heston  characteristic function hides  the analycity
property and difficults its utilization. So  we first give a new expression of the characteristic function, but more
importantly,
we use that the  analycity property is equivalent that the random variable has (real) moment generating function in a neighborhood
of the origin: there is $\varepsilon>0$ such that
$$M(u)=\E[e^{uX_t}]<\infty,\  -\varepsilon <u<\varepsilon.$$
In that case, the function $M(u)$ is a  real analytic function in $(-\epsilon,\epsilon)$
 and the main properties of the characteristic function can be studied through  $M(u)$, which is
a real function, and in general much  simpler to analyze.

As a consequence of the study of   the moment generating function we obtain the domain of that function and
we give a numerical simple procedure to compute the poles of the characteristic function in its strip of convergence.
This has several practical consequences and  we apply
it to the computation of the smile wings parameters in the
formula
given by Lee \cite{Lee04}. We also apply these results to the assessment of the moments of a particular
model;    such study is complementary of the one of Andersen and Piterbarg \cite{AndPie07}.

The second contribution of this paper  is that we factorize the moment generating function of   log-spot.
This allows us to identify a Bessel type densities as the building blocks of the log-spot.
This result  has   interesting applications. For example, it allows us to construct
a sequence of random variables converging in law to  log-spot.
On the other hand, for a certain combinations of the parameters,  log-spot $X_t$  is a sum
 of independent non-centered $\chi^2$ random
variables, and we can identify $X_t$   as a member of the non-homogeneous second Wiener chaos (see Janson
 \cite[chapter 6]{Jan97}); this agrees with the intuition that for some parameters, the stochastic volatitlity, given by
 a CIR model (Cox {\it et  al.} \cite{CoxIngRos85}), has the law of the sum of the squares of a finite number of
  independent Ornstein-Uhlenbeck processes, that are in the second Wiener chaos, and such property is transferred
  to $X_t$.

  The paper is organized as follows. First we deduce the moment generating function and  the characteristic
  function of log-spot; as far as we know, these expressions are new. Both can be obtained manipulating the expressions
  obtained by Dufresne \cite{Duf08} or the standard  expressions of the characteristic functions
  (see for example, Gatheral \cite{Gat05} or  Albrecher {\it et al.} \cite{AlbHanMaySch07}), however we prefer
  to give a new deduction
  from scratch since our procedure is quite general and can be applied to other problems. In Section \ref{domini},
  we obtain  the domain of the moment
  generating function. In Section \ref{applications}, we give some applications.
Finally, in Section \ref{factorsec},
   using techniques of complex analysis, the  moment generating
  function of log-spot is factorized. The analysis of the factors   allows us  to identify them
    as  moment generating functions of  Bessel type random variables, and  to construct a sequence of random variables
    that converges in law to {\it log} spot. In the Appendix we review some facts on moment generating function of
    a random  variables and we put technical details of the  proofs.

\section{The Heston model}
The Heston model \cite{Hes93} is defined by the system of stochastic differential equations
 \begin{equation}
 \label{heston}
 \left.\begin{aligned}
dS_t&=\mu S_t\, dt+S_t\sqrt{V_t}\, dZ(t)\\
dV_t&=a(b-V_t)\,  dt+c\sqrt{V_t}\, dW(t)
\end{aligned} \right\}
\end{equation}
with initial conditions $S_0=s_0>0$ and $V_0=v_0\ge 0$, where $a,\, b>0$ and  $c\in \R-\{0\}$ are constants, and $W$
  and $Z$ are two standard correlated Brownian
motions,
$\langle Z,W\rangle_t=\rho t,$ for some $\rho\in[-1,1].$
 The process $V_t$ is a Feller diffusion (Feller \cite{Fel51}) or, in
the financial literature,
 a CIR model (Cox {\it et al.} \cite{CoxIngRos85}). The parameter
$a$ is called the {\it the mean reversion factor},  $b$ is called the {\it long term volatility} and it
 is also written $V_\infty$, and  $c$ is called the {\it vol-of-vol}. Write
$$X_t=\log S_t-\mu\, t.$$
By the It\^{o} formula,
$$dX_t=-\frac{1}{2} V_t\, dt+\sqrt{V_t}\, dZ_t,$$
with initial condition $X_0=x_0=\ln s_0$. Thus, we will consider the system  \begin{equation}
 \label{heston-log}
 \left.\begin{aligned}
dX_t&=-\frac{1}{2} V_t\, dt+\sqrt{V_t}\, dZ_t\\
dV_t&=a(b-V_t)\,  dt+c\sqrt{V_t}\, dW(t)
\end{aligned} \right\}
\end{equation}

\subsection{The moment generating function of $\bs{X_t}$}
First, we check that $X_t$ (indeed $(X_t,V_t)$) has moment generating function, and later  we deduce its expression as the solution
of a (real)
PDE obtained from It\^{o} formula.

 For every $u,v\in \R$, the random variable $e^{uX_t+vV_t}$  is positive, so we can compute its expectation but it can be
infinite.
Write
$$Z_t=\rho W(t)+\sqrt{1-\rho^2}W'(t),$$
where $W'$ is a Brownian motion independent of $W$, and use the habitual trick
$$g(u,v,t):=\E\big[e^{uX_t+vV_t}\big]=\E\Big[\E\big[e^{uX_t+vV_t}/W'(s),\ 0\le s\le t\big]\Big].$$
  We obtain
\begin{align*}
g(u,v,t)&=
\exp\big\{x_0 u-\frac{u v_0\rho }{c}-\frac{u\rho ab t}{c}\big\}\\
&\cdot \E\Big[\exp\Big\{\big(v+\frac{\rho u}{c}\big)V_t+\big(\frac{u^2}{2}-\frac{u}{2}-\frac{\rho^2u^2}{2}+\frac{u a\rho}{c} \big)
\int_0^t V_s\, ds\Big\}\Big].
\end{align*}
Note that the coefficient of $\int_0^t V_s\, ds$ is the equation of a parabola in $u$ through  the origin, so if $u$ is near zero, that coefficient
should be also near zero. Since both $V_t$ and $\int_0^t V_s\, ds$ have moment generation function (see Dufresne
\cite{Duf01})
it follows that the expectation is finite for $(u,v)$ in a neighborhood of $(0,0)$.
Fix $T>0$.
Applying  the It\^{o} formula  to $\exp\{uX_t+vV_t\}$ and taking expectations, we have
\begin{align*}
g(u,v,t)&=e^{ux_0+vv_0}+v a b \int_0^t g(u,v,s)\, ds
+\Big(- \frac{u}{2}
+\frac{u^2}{2} +
\frac{v^2c^2}{2}
+uv\rho c-va\Big)\int_0^t \frac{ \partial g(u,v,s)}{\partial v}\, ds,
\end{align*}
where we have used the property of the moment generating function
$$\frac{ \partial g(u,v,s)}{\partial v}=\E\big[e^{uX_s+vV_s}V_s\big]$$
Differentiating with respect to $t$,   we get
$$\frac{\partial g(u,v,s)}{\partial t}-p(u,v)\frac{\partial g(u,v,t)}{\partial v}=a b v g(u,v,t),$$
where
$$p(u,v)=- \frac{u}{2}+\frac{u^2}{2} + \frac{v^2 c^2}{2} +uv\rho c-v a.$$
This equation has a unique solution that is
$$g(u,v,t)=\Bigg(\frac{p\big(u,\phi(u,v,t)\big)}{p(u,v)}\Bigg)^{ab/c^2}\exp\big\{ux_0+\phi(u,v,t)v_0-
\frac{ ab t\rho u}{c}+\frac{a^2b t}{c^2}\big\},$$
where
\begin{align*}
P(u)&=\sqrt{(a-\rho c u)^2+ c^2(u-u^2)},\\
\gamma(u,v)&=-2\arctan \!{\rm h}\big((c^2v+c\rho u-a)/P(u)\big)/P(u)\\
\varphi(u,v,t)&=-\frac{\rho u}{c} +\frac{a}{c^2}-\frac{1}{c^2}P(u)\tanh\Big(P(u)\big(t+\gamma(u,v)\big)/2\Big).
\end{align*}
For $v=0$, we get the moment generating function of $X_t$, $M_t(u)$; when there is no confusion  we will
 suppress the subindex $t$
and write $M(u)$. After  some tedious manipulations,
$M(u)$ can be written as
 \begin{align}
 \label{generating}
M(u)&=\E\Big[\exp\big\{ uX_t\big\}\Big]\notag \\
&= \exp\{x_0u\}
\ \bigg(\frac{e^{(a-c\rho u) t/2}}{\cosh(P(u)t/2)+\big(a -c\rho u)\sinh(P(u)t/2)/P(u)}\bigg)^{2a b/c^2}\notag\\
\notag\\
&\qquad
\cdot\exp\bigg\{-v_0\,\frac{(u-u^2)\sinh(P(u)t/2)/P(u)}
{\cosh(P(u)t/2)+\big(a-c\rho u \big)\sinh(P(u)t/2)/P(u)}\bigg\},
\end{align}
where
$$P(u)=\sqrt{(a-\rho c u)^2+ c^2(u-u^2)}.$$

\begin{remark} Formula (\ref{generating}) coincides with the one that can be deduced from the joint
Laplace-Mellin transformation of $S_t,V_t,\int_0^tV_s\, ds$ given by Dufresne \cite[Theorem 12]{Duf08}.
\end{remark}

\subsection{The characteristic function of $\bs{X_t}$}
\label{complexm}

For $z$ complex, consider the function
\begin{align}
 \label{laplace}
\Phi(z)=&
\exp\{x_0 z\}
 \ \bigg(\frac{e^{(a-c\rho z)t/2}}{\cosh(P(z)t/2)+(a -c\rho z)\sinh(P(z)t/2)/P(z)}\bigg)^{2a b/c^2}\notag \\
\notag \\
& \qquad \cdot\exp\bigg\{-v_0\ \frac{(z-z^2)\sinh(P(z)t/2)/P(z)}
{\cosh(P(z)t/2)+(a-c\rho z)\sinh(P(z)t/2)/P(z)}\bigg\}.
\end{align}
Write $p(z)=(a-\rho c z)^2+ c^2(z-z^2)$ the second degree polynomial within $P(z)$. Since $p(0)=a^2>0$, using
standard techniques of complex analysis, we see that $\Phi(z)$ is well defined and analytic in a neighborhood of $0$.
Obviously, $\Phi(u)=M(u)$ on a (real) neighborhood of 0. Then (see the Appendix, Proposition A.\ref{coinc}, and
Section \ref{factorsec}),
  the characteristic function of $X_t$   is $\Phi(iu)$.
Explicitly, for $u\in \R$,
 \begin{align}
 \label{caract1}
 \varphi(u)=\E[e^{iuX_t}]=\Phi(iu)&=
 \exp\{ix_0u\}
\ \bigg(\frac{e^{a t/2} }{\cosh(dt/2)+\xi\sinh(dt/2)/d}\bigg)^{2ab/c^2}\notag \\
&\qquad  \cdot \exp\bigg\{-v_0\,\frac{(iu+u^2)\sinh(dt/2)/d}
{\cosh(dt/2)+\xi\sinh(dt/2)/d}\bigg\},
\end{align}
where
\begin{align*}
d&=d(u)=P(iu)=\sqrt{(a-c\rho iu)^2+c^2(iu+u^2)},\\
\xi & =\xi(u)=a-c\rho ui.
\end{align*}
After some computations we arrive to the formula of Albrecher {\it et al.} \cite{AlbHanMaySch07}
 \begin{equation}
 \label{caract2}
 \varphi(u)=\exp\{i x_0 u\}
 \exp\Big\{\frac{ab}{c^2}\Big( (\xi-d)t-2 \log\frac{1-ge^{-dt}}{1-g}\Big)\Big\}
\, \exp\Big\{\frac{v_0}{c^2}(\xi-d)\frac{1-e^{-dt}}{1-ge^{-dt}}\Big\},
\end{equation}
where $$g=g(u)=\frac{\xi-d}{\xi+d}.$$
Of course, formula (\ref{caract1}) looks more complex than the compact (\ref{caract2}). However, when one recovers
from the shock, one  realizes that the former is easier to handle that the latter.

\subsection{Inversion of Heston model and some comments on  the parameters}
We will say that  a  process $S=\{S_t,\ t\in [0,T]\}$ is a Heston type process, and  we write  $S\sim {\cal HP}_{\Prob}(a,b,c,\rho,s_0,v_0,\mu)$ if $ S_t$ verifies a system (\ref{heston-log}) with some stochastic
volatility $V_t$. When there is no confusion with the underlying probability $\Prob$ we will omit it.

 Recent results of del Ba\~{n}o \cite{Ban08} show that if $S\sim {\cal H P}_{\Prob}(a,b,c,\rho,s_0,v_0,\mu)$, $a>c\rho$, then
\begin{equation}
\label{inversion}
S^{-1}\sim {\cal H}_{\Q}(a-c\rho  ,ab/(a-c\rho),c,-\rho,s_0^{-1},v_0,-\mu),
\end{equation} where $\Q$ is the probability given by
$$\frac{d \Q}{d\Prob}=e^{X_T}.$$
To prove that property it is needed to work with the whole process $S$. However, an easy verification
can be done  using the moment generating function (\ref{generating}).
We will use this property in Section \ref{domini} to compute the values of $\E[\exp\{uX_t\}]$ for negative values
of $u$.

It can  also be proved that
$${\cal H P}_{\Prob}(a,b,c,\rho,s_0,v_0,\mu) \sim {\cal H P}_{\Prob}(a,b,-c,-\rho,s_0,v_0,\mu).$$
Again, it is needed to consider  the  process to prove this equality, but a check is deduced  from the
 expression of
the moment generating function (\ref{generating}).
So, without loss of generality  we will assume from now on that $c>0$.

\section{The domain of  the moment generating function}
\label{domini}
In this section we deduce the domain of the moment generating function (\ref{generating}); this deduction is
 not direct since we obtained  $M(u)$ not by the computation of  the expectation $\E[e^{uX_t}]$ but by an indirect way.
 So, we only know that the moment generating function coincides with the function given in the right hand side
 of (\ref{generating}) in a neighborhood of zero.  However, as stated by Lemma 3.3 below,
since the function is analytic, we are in safe land. In the first subsection, using this idea we do a first
 study of the moment generating
 function. Later, in the second subsection we work with the analytic continuation of the function  in the  right
 hand side of  (\ref{generating}).
\subsection{Preliminary  study of   the domain of the moment generating function}
\label{primer-estudi}

The main ingredient  of the moment generating function given in (\ref{generating}) is
\begin{equation}
\label{discriminant}
f(u):=\cosh(P(u)t/2)+(a-c\rho u)\,\frac{\sinh(P(u)t/2)}{P(u)},
\end{equation}
where
$$P(u)=\sqrt{(a-\rho c u)^2+ c^2(u-u^2)}.$$
Write
\begin{equation}
\label{segongrau}
p(u):=(a-\rho c u)^2+ c^2(u-u^2).
\end{equation}
When $\rho\ne \pm1$, $p(u)$ represents a parabola
with leading coefficient $c^2(\rho^2-1)\le 0$, and for $u=0$, we have $a^2>0$.
So,  it is an inverted parabola with real roots,
$u_-<0<u_+$, given by
\begin{equation}
u_{\pm}=\frac{c-2a\rho\pm\sqrt{4a^2+c^2-4ac\rho}}{2c(1-\rho^2)}.
\label{umes}
\end{equation}
When $\rho=-1$, then  $$p(u)=a^2+c(c+2a)u,$$ is a straight line with positive slope
that intersects the horizontal axis at $u_-=-a^2/(c(c+2a))<0,$
and we write $u_+=\infty$.

Similarly, for $\rho=1$, $p(u)$ degenerates in the straight line
$$p(u)=a^2+c(c-2a)u,$$
and the slope can be negative, postitive or  zero or negative, and
\begin{itemize}
\item When $2a>c$, then $u_-=-\infty$.
\item  When $2a<c$, then $u_+=\infty$.
\item When $2a=c$, then $u_-=-\infty$ and $u_+=\infty$.
\end{itemize}

Thus, for every $\rho\in[-1,1]$, we have that $p(u)\ge 0$ on $[u_-,u_+]$, hence the function $f(u)$ is well defined
and analytic in such (possible infinite) interval.
Denote by $D(X_t)$ the domain of $M(u)$,
 $${D}(X_t)=\big\{u\in \R: \ M(u)=\E\big[ e^{uX_t}\big]<\infty\big\}.$$
In principle (see Lemma \ref{extensio} bellow)  $M(u)$ is defined in the subinterval of
$[u_-,u_+]$
 between  the biggest   negative zero and the smallest positive zero of $f(u)$.
Next proposition summarizes the study of such zeroes in that interval.

\begin{proposition}
\label{primerafita}
For every $t>0$, and $\rho\in[-1,1]$ there is the inclusion $[u_-,1]\subset D(X_t)$. Moreover,
\begin{enumerate}[\bf 1.]
\item When $a\ge c\rho$ (in particular, for every $\rho<0$), then for all $t>0$, the function $f(u)$
has no zeroes in $[u_-,u_+]$, and consequently,  $[u_-,u_+]\subset D(X_t).$ (Except for $\rho=1$ and $c=2a$).
\item When $a <c\rho$, write  $t_0=2/(c\rho u_+-a)\ge 0$
\begin{enumerate}[(i)]
\item If $t<t_0$, then $[u_-,u_+]\subset D(X_t)$, and  $f(u)$
has no zeroes in $[u_-,u_+]$.
\item If $t\ge t_0$, then  $f(u)$
has no zeroes in $[u_-,1]$, and has one and only one zero in $(1,u_+]$.
\end{enumerate}
\item When $\rho=1$ and $a=2c$, then $D(X_t)=\big(-\infty,1/(1-e^{-at})\big)$.
\end{enumerate}

\end{proposition}

\begin{remark}

\rule[0mm]{0cm}{1cm}

\begin{enumerate}[\bf 1.]
\item The cases $\rho=\pm 1$ are   specially important. We stress that
\begin{enumerate}[(i)]
\item For $\rho=-1$, we are always in case
$1$ and $u_+=\infty$.
\item  For $\rho=1$,   when $ a\ge c$, we are in case 1 (except if $a=2c$). For
 $a<c$ we are  in case 2; however, if $a<c<2a$, then $u_+=\infty$, and  we have
 $t_0=0$, and thus we are in case (ii) for all $t>0$.
\end{enumerate}

 \item From the preceding proposition it follows that $\E[S_t]<\infty$
for every $a,b,c>0$ and $\rho\in[-1,1]$.
 Moreover, when $\mu=0$,  by construction, $\{S_t, \, t \in [0,T]\}$ is an exponential local martingale,
that is a positive supermartingale, see Revuz and Yor \cite[pages 148 and 149]{RevYor99}, and
$$\E[S_t]=M(1)=e^{x_0},$$
so $\{S_t, \, t \in[0,T]\}$ is a true martingale.  This was   proved by Andersen and Piterbarg \cite[Proposition 2.5]{AndPie07}
 using the Feller explosion criteria and Girsanov Theorem.

\item As a continuation of the  preceding point,  we should remark that  $\{S_t, \, t \in[0,T]\}$  is not always a square
 integrable martingale.
This can cause some problems. See Subsection \ref{segonordre}.
\end{enumerate}
\end{remark}
\bigskip

\subsection{Computation of the absciss\ae\ of convergence  of the moment generating function}
\label{extension}
When there is no root of $f(u)$ in $[1,u_+]$ the domain $D(X_t)$ is larger than $[u_-,u_+]$ (in particular, when $a>\rho c$).
To carry out this study,  we need more properties   of the moment generating function. Consider an arbitrary
 random variable $X$, with   moment generating function $M_X(u)$,
and   domain $D_X=\big\{u\in \R: \ M_X(u)<\infty\big\}$.
Remember that $D_X$
 is a interval of $\R$ (finite or infinite, open or closed from one side or the other, that always include the origin,
and it may be just $\{0\}$), and
 $M_X$ is analytic in the interior of ${D_X}$.
 The left (respectively right)
 extreme of ${ D_X}$ is called the left (resp. right) abscissa of convergence, and plays a major role.
The following property is well known, but since
 it is key in this paper, we stress it. We give the property for the right abscissa of convergence, and
 a similar statement is true for the left abscissa,

 \begin{lemma}
\label{extensio}
  Let $X$ a random variable such that there is a neighborhood of zero included in $D_X$. Assume that $(r,s] \subset { D_X}$, and that there is an analytic function $h:(p,q)\to \R$ such that
 \begin{enumerate}
 \item $(r,s]\subset (p,q)$.

 \item $M_X(u)=h(u)$,\  \text{for}\  $ u\in (r,s]$.

 \end{enumerate}
Then $M_X=h$ on $(p,q)$. Moreover, if \, $\lim_{u\nearrow  q} h(u)=\infty,$ then
the right--abscissa of convergence of $M_X$ is the point $q$.

 \end{lemma}

\noindent{\it Proof.}

\medskip

 Denote the interior of $D_X$ by $(\alpha,\beta)$. If $\beta=\infty$,  by analytic continuation, $M_X=h$ on $(r,q)$.
 Consider $\beta<\infty$, so
  $\beta$ is the finite right--abscissa of convergence.  Then the function  $M_X$
  has a singularity at  $\beta$ (see the Appendix, Proposition A.\ref{abscisa})
 Thus,    $s\ne \beta$, because $M_X(s)=h(s)$ and $h$ is analytic in $s$. Hence, $\beta>s$.
 In the same way, $\beta<q$ is contradictory,
then $\beta\ge q$ and by analytic continuation, $M_X=h$ on $(r,q)$.
The second part of the Lemma is obvious.
$\qquad \square$

\bigskip

Now we apply the preceding lemma to the moment generating function  of the log-spot, $M(u)$.
When the function $f(u)$ given in (\ref{discriminant}) has no zeros in $[1,u_+]$,  since
$$\lim_{u\searrow u_-} M(u)<\infty \qquad \text{and}\qquad \lim_{ u\nearrow u_+} M(u)<\infty,$$
 by Lemma \ref{extensio}, the domain of $M(u)$, $D(X_t)$, is bigger than $[u_-,u_+]$
Then, to assess more carefully that domain,  consider the function $\cosh\sqrt x$, for $x>0$. Its Taylor expansion is
$$\cosh\sqrt x=\sum_{n=0}^\infty \frac{x^n}{(2n)!}.$$
The series on the right defines an entire function, say $L_1(x)$. However, when $x<0$, that  series
coincides with the Taylor expansion
of $\cos \sqrt{- x }$. Hence,  $L_1(x)$ is an entire function that, when written   as the composition of elementary functions,
has different expression according to whether $x> 0$ or $x< 0$, that is
$$L_1(x)=\begin{cases}\cosh \sqrt{x}, & \text{if $x\ge 0$},\\
\\
\cos \sqrt{- x}, & \text{if $x\le 0$}.
\end{cases}$$

In a similar way, the function $(\sinh \sqrt x)/\sqrt x,\ x>0,$
 can be analytically continuated  to   negative numbers as $(\sin \sqrt{- x})/\sqrt {-x},\ x<0$.

\bigskip

Denote by $u_-^*\ge -\infty$  the left abscissa of convergence of $M(u)$  and by $u_+^*\le \infty$ the
right abscissa.  Put
$$\wt P(u)=\sqrt{-p(u)}.$$ For $u\in (u_-^*, u_-)$  or $u\in (u_+,u_+^*)$,
by Lemma \ref{extensio}
the moment generating function is
 \begin{align}
 \label{generating2}
M(u)&=
\exp\{x_0u\}\
 \bigg(\frac{e^{(a-c\rho u) t/2}}{\cos(\wt P(u)t/2)+\big(a -c\rho u)\sin(\wt P(u)t/2)/\wt P(u)}\bigg)^{2a b/c^2}
  \notag\\
\notag\\
&\qquad \cdot\exp\bigg\{-v_0\,\frac{(u-u^2)\sin(\wt P(u)t/2)/\wt P(u)}
{\cos(\wt P(u)t/2)+\big(a-c\rho u \big)\sin(\wt P(u)t/2)/\wt P(u)}\bigg\}.
\end{align}

\begin{figure}[htb]
\centering
%
%
%
%
%

\includegraphics[scale=0.8]{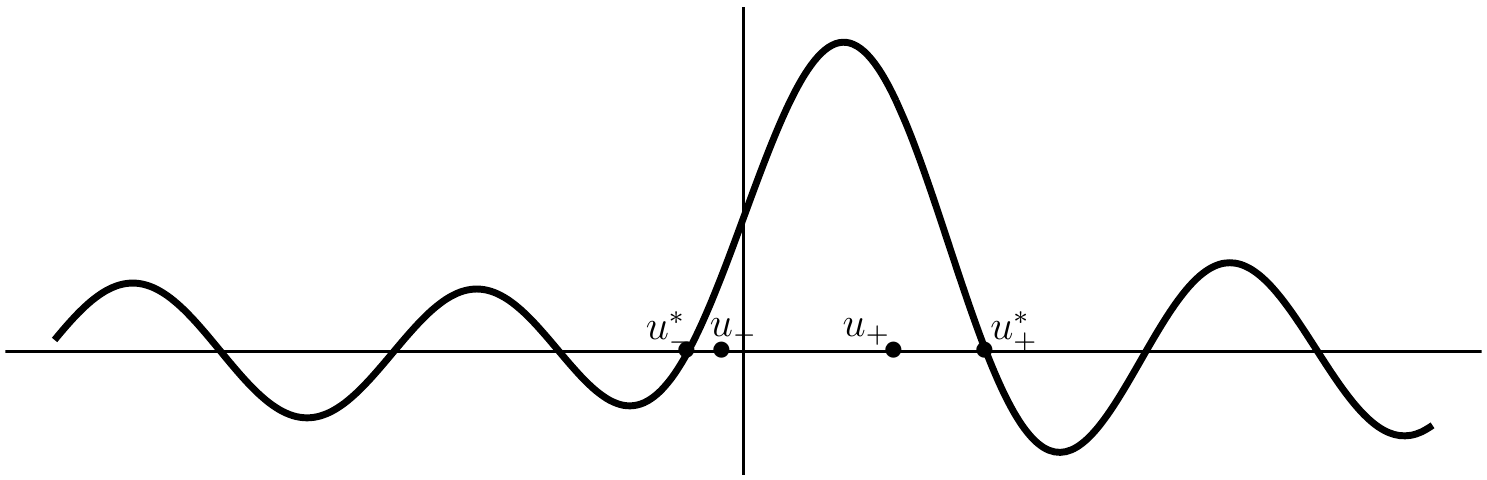}

\caption{Plot of $F(u)$.
 The points  $u_-^*$ and $u_+^*$ are the abciss\ae\  of convergence of $M(u)$.}
\label{F(u)}
\end{figure}

 In that expression, the main part is the function
\begin{equation}
\label{funcioq}
 \wt f(u):=\cos(\wt P (u)t/2)+(a -c\rho u)\frac{\sin(\wt P(u)t/2)}{\wt P(u)}.
 \end{equation}
 Both $f$ and $\wt f$ are defined in disjoint sets, and can be combined in a new function
\begin{equation}
\label{unificacio}
F(u)=
\begin{cases} f(u), & \text{if $u\in[u_-,u_+]$,}\\
\\
\wt f(u),  & \text{if $u < u_-$ or $u>u_+$,}
\end{cases}
\end{equation}
that is analytic in $\R$.
See Figure \ref{F(u)} for a plot of that function.

\begin{figure}[htb]
\centering

\includegraphics[scale=0.8]{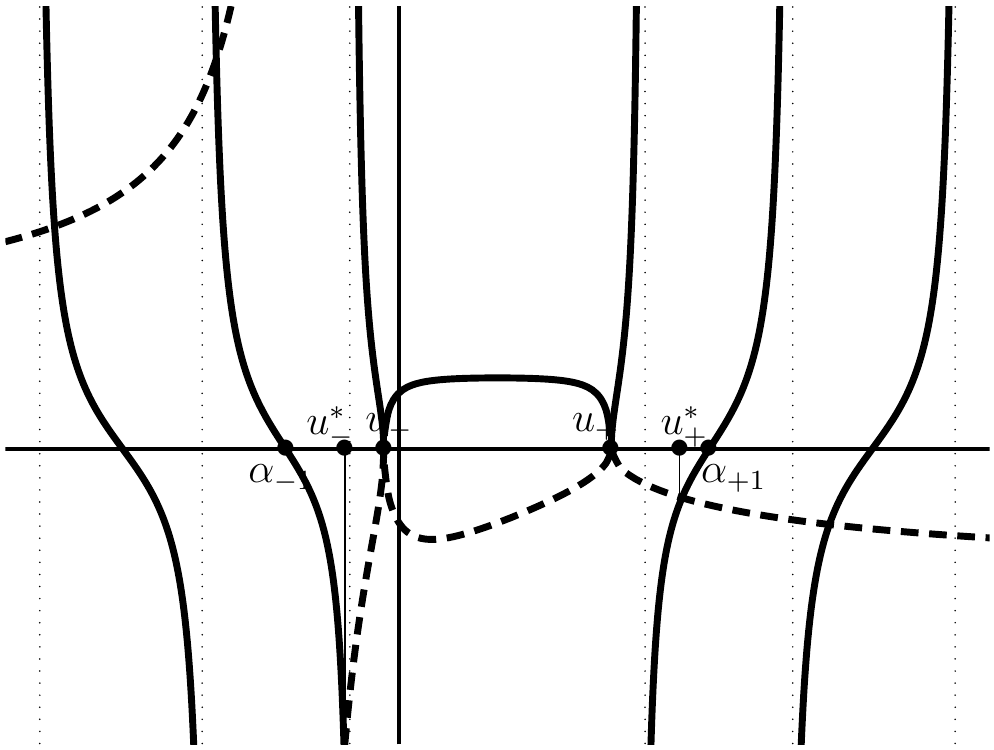}
\caption{Solid line: Plot of $\tanh\big(P(u)t/2\big)$, for $u\in [u_-,u_+]$ or $\tan\big(\wt P(u)t/2\big)$ in the complementary.
Dashed line: Plot of $-P(u)/(a-\rho c u)$ or $-\wt P(u)/(a-\rho cu)$.}
\label{interseccions}
\end{figure}

Then, to find the right abscissa of convergence, we need to find the  zero of the function
$f(u)$ in $[1,u_+]$, and if there is no zero, then  to find the smallest  zero $u>u_+$
  of $\wt f(u)$. Note also that the real zeroes of $\wt f(u)$ are the real solutions of the equation (see Figure
  \ref{interseccions})
\begin{equation}
\label{equacioq}
 \tan(\wt P(u)t/2)=-\frac{\wt P(u)}{a-c\rho u}.
 \end{equation}
 (Except   when $a/c\rho=k\pi/2$ for some natural number $k$).

 For the left abscissa of convergence we need to look for the
 biggest  solution $u<u_-$ of $\wt f(u)=0$, or, equivalently, to work with equation (\ref{equacioq}).

In order to give a  bound for the abscisses of convergence, denote by $\alpha_{\pm 1}$ the solutions of the equation
$$p(u)=-\frac{4\pi^2}{t^2},$$
that for all $t>0$  are real and $\alpha_{-1}<u_-$  and $\alpha_{+1}>u_+$. Also put  $\beta_{\pm 1}$ the solutions
of  $p(u)=-2\pi^2/t^2.$ Consider the following cases:
\begin{enumerate}[\bf 1.]

\item If $a>\rho c$, then $u_+^*\in(u_-,\alpha_{+1})$. This can be deduced from the consideration that
the image of the function $\tan(\wt P(u)t/2)$ on that interval is $\R$, and the properties of the   function
in the right hand side of
(\ref{funcioq}). The only  case not clear is when $a/(c\rho)=\pi/2$, due to the fact that
$-{\wt P(u)}/(a-c\rho u)$ has a vertical asymptote at $u=\beta_1$; this case is studied by direct
inspection.

\item If $a<\rho  c$, then remember that $f(1)>1$, and, on the other hand,
$$\wt f(\beta_{+1})=(a-c\rho \beta_{+1})2/\pi<0,$$
because $\beta_{+1}>1$ and  $\rho>0$. So $F(u)$ has at least one root in $(1,\beta_{+1})$.

\end{enumerate}

Joining these comments with Proposition \ref{primerafita} we have

\begin{theorem}
\label{segonafita}
With the above notations,
\begin{enumerate}[\bf 1.]
\item If $a>\rho c$,  the right abscissa of convergence  $u_+^*$ is  the smallest  zero $u>u_+$ of $\wt f(u)=0$, and
  $u_+^*\in(u_+,\alpha_{+1})$. (Except for $\rho=1$ and $a=2c$).
  \item If $a=\rho c$, the right abscissa of convergence is $u_+^*=1$.

 \item If $a<\rho  c$,  let  $t_0=2/(c\rho u_+-a)\ge 0$.
\begin{enumerate}[(i)]
\item If $t<t_0$, then $u_+^*$ is smallest  zero $u>u_+$ of $\wt f(u)=0$, and $u_+^*\in(u_+,\beta_{+1})$.
\item If $t\ge t_0$, then $u_+^*$ is the zero of $f(u)$ in $(1,u_+]$.

\end{enumerate}
\item In every case, the left abscissa of convergence, $u_-^*$ is the biggest zero $u<u_-$ of $\wt f(u)=0$
and $u_-^*\in(\alpha_{-1},u_-)$.
\end{enumerate}
\end{theorem}

\begin{remark}

\rule[0mm]{0cm}{1cm}

\begin{enumerate}

\item When $a>\rho c$ the inversion formula (\ref{inversion}) can be used to compute $u^*_-$ in terms of
$u^*_+$ of the inverted model. Specifically,
$$u^*_-(a,c,\rho,t)=-u^*_+(a-c\rho,c,-\rho,t)+1.$$

\item This theorem gives a  direct procedure  to invert the formulas  of   Andersen and Pitebarg \cite[Proposition 3.1]{AndPie07}.

\end{enumerate}
\end{remark}

\section{Applications}
\label{applications}

In this section we present some applications of the exact knowledge  of the abciss\ae\  of convergence of the
moment generating function of log-spot.

\subsection{The smile at extreme strikes}
One of the motivations for the study of the abciss\ae\   of convergence of the
moment generating function
 of log-spot in the Heston model in this work is the outstanding result of Roger Lee (\cite{Lee04}) where an
explicit relation is found between the assymptotic behaviour of the volatility smile and these  abciss\ae.
This can be of interest in designing sensible smile interpolation and extrapolation schemes as has been
 shown by  Gatheral  \cite{Gat04}.
In the case under consideration,
that of the Heston dynamics, if $u_+^*$ and $u^*_-$ are the abciss\ae\  of convergence of log-spot,
 then according to Lee \cite{Lee04} the asymptotic behaviour of the Heston smile for expiry $T$
 as the strike $K$ goes to infinity is $\sigma(K)\approx \sqrt{\beta_R \ln(K)/{T}}$ where $\beta_R\in[0,2]$
  is defined by
$$
\frac{1}{2\beta_R}+\frac{\beta_R}{8}-\frac{1}{2}=u_+^*-1.
$$
Likewise the behaviour as $K$ approaches zero is $\sigma(K)\approx \sqrt{-\beta_L  \ln(K)/{T}}$ where $\beta_L\in[0,2]$ is defined by
$$
\frac{1}{2\beta_L}+\frac{\beta_L}{8}-\frac{1}{2}=-u^*_-,
$$
As an example, consider  market data for a one year equity smile   $\rho=-90\%, \, a=2,\, c=80\%, \, b=15\%^2$ and $t=1$.
The calculations outlined above yield $u_+^*=37.43$ and $u_-^*=-3.21$. By Roger Lee's formulas above this implies
the behaviour of the smile wings is ruled by
$$\beta_R=0.01\quad\text{and}\quad \beta_L=0.13$$
where, as expected, the smile at low strikes has a larger coefficient.
These numbers can be of use to design an extrapolation scheme for the Heston smile in extreme strikes
where the numerical integration breaks down.

\subsection{The importance of second order moment of spot in  stochastic volatility models}
\label{segonordre}
The moments of spot and the abciss\ae\  of convergence of log-spot are related because
\begin{equation}
\label{sup-inf}
u_+^*=\sup\{u\in \mathbb{R}:\ \E[S_t^u]<\infty\}\quad\text{and}\quad u_-^*=\inf\{u\in \mathbb{R}:\ \E[S_t^{u}]<\infty\}.
\end{equation}
So, if $u_+^*<2$, then
the second order moment of spot is infinite.
This can cause  problems in pricing certain standard
European derivatives,  as the next couple of examples show.

\bigskip

\noindent{\bf  1. Pricing an FX performance note.}
Consider a performance note that pays Euros on the performance of the EUR/USD exchange rate.
This is a contract that pays the following amount in Euros at expiry $T$
$$
\text{Payoff}(EURUSD_T)=\text{Notional}_{EUR} \left( \frac{EURUSD_T}{EURUSD_0} -1 \right) \qquad EUR
$$
The notation should be self explanatory: $EURUSD_t$ is the EUR/USD exchange rate at time $t$ and  $t=0$ is today.
By the standard martingale methods the price of a derivative product is simply its discounted expectation,
this might (and does) lead naive market participants to price such a transaction as
$$
\text{Present Value} = Df_T^{EUR} \text{Notional}_{EUR} \left( \frac{F_T}{EURUSD_0} -1 \right) \qquad EUR
$$
where we have used the fact that the risk neutral expectation of spot is the forward $F_T$, and write $Df_T^{EUR}$ for the relevant discount factor.
Unfortunately this approach is flawed since EUR/USD is the price in US Dollars of one Euro and in the equations above we are basing a Euro payment on this quantity.
Its correct price can be derived by noting that its payout in Dollars is
$$
\text{Payoff}(EURUSD_T)=\text{Notional}_{EUR} \left( \frac{EURUSD_T}{EURUSD_0} -1 \right) \cdot EURUSD_T \qquad USD
$$
and here we can apply the martingale methods to yield
\begin{eqnarray*}
\text{Present Value} &=&Df_T^{USD} \text{Notional}_{EUR} \mathbb{E} \left( \left(
\frac{EURUSD_T}{EURUSD_0} -1 \right) \cdot EURUSD_T \right)\\
&=&Df_T^{USD} \text{Notional}_{EUR}  \left( \frac{\mathbb{E}\left(EURUSD_T^2\right)}{EURUSD_0} -F_T \right) \\
\end{eqnarray*}
an expression that involves the second moment of spot.
A model that has an infinite second moment for spot (and these things do appear in practice)
 will price such a contract at infinity.

\bigskip

\noindent{\bf  2. LIBOR paid in arrears.}
A similar type of deal occurs in the fixed income derivatives market under the name of LIBOR paid in arrears.
Normally LIBOR is a rate fixed on a certain date $T_1$ and paying at a later date $T_2$.
A simple contract depending on LIBOR is a FRA (Forward Rate Agreement) whose payout is defined by
$$
\delta ( L-K )\qquad \text{paid at time }T_2
$$
here L is the LIBOR rate which fixes at $T_1<T_2$ and $\delta $ is the day-count fraction which is approximately $T_2-T_1$.
If we use the $T_2$-forward measure, the value of such a contract is simply the discounted expectation $\delta P(0,T_2)\mathbb{E}(L-K)$ where $P(0,T_2)$ is the $T_2$-discount bond.
Given that $L_t=\frac{1}{\delta}(P(t,T_1)-P(t,T_2))/P(t,T_2)$ is the ratio of a tradeable instrument by the numeraire, we know it is a martingale in the $T_2$-forward measure and since $L_{T_1}$ is simply the LIBOR rate, the price of our FRA is
\begin{eqnarray}\label{FRA}
P(0,T_2)\delta \mathbb{E}( L-K )&=&
P(0,T_2)\delta \mathbb{E}( L_{T_1}-K )\\
&=&
P(0,T_2)\delta ( L_0-K )\\ \nonumber
\end{eqnarray}
where $L_0=\frac{1}{\delta}(P(0,T_1)-P(0,T_2))/P(0,T_2)$ is the so-called forward LIBOR rate,
the rate that makes the FRA have zero value.
A LIBOR in-arrears transaction is based on LIBOR paid \emph{at the wrong time}, its payout being
$$
\delta ( L-K )\qquad \text{paid at time }T_1
$$
A naive idea to price these transactions is to simply discount the by the forward discount factor
$P(0,T_2)/P(0,T_1)$ which of course just affects equation
(\ref{FRA}) by a multiplicative factor and does not alter the fair forward price $L_0$.
This is wrong because the expectation of $L$ is no longer the forward LIBOR in the $T_1$-forward measure which this approach implicitly uses. Some banks have
been arbitraged in the past by using this naive approach.
As in the case of the EURUSD note described above, if
we convert the payment to a payment at time $T_2$ then we can use the expectations in the $T_2$-forward measure, the payout at time $T_2$ is simply the accrued amount
$$
\delta (1+\delta L)\cdot ( L-K )\qquad \text{paid at time }T_2
$$
and its price will be
$$
P(0,T_2)\delta \mathbb{E}\left((1+\delta L)\cdot ( L-K )\right)
$$
and so the fair strike is strictly larger than the FRA rate
$$
L_0+\delta \frac{\mathbb{V}ar(L)}{1+\delta L_0} $$
an expression that involves the second moment of the LIBOR rate.

\subsection{Dependence of the absciss\ae\ of convergence on the time.}
In this subsection we assume that $\mu=0$, so $\{S_t,\, t\in[0,T]\}$ is  a  martingale.
In order to study the dependence of the abciss\ae\  of convergence on the time, we denote by $u_+^*(t)$ and
$u_-^*(t)$ the abciss\ae\  for $X_t$. The following proposition is a general property of a positive martingale.

\begin{proposition}
Consider $0\le t<t'\le T$. Then
$$u_+^*(t)\ge u_+^*(t') \quad\text{and}\quad u_-^*(t)\le u_-^*(t').$$
\end{proposition}

\noindent{\it Proof.}

Fix  for a moment $r\ge 1$. The function $\phi(x)=x^r$ on $(0,\infty)$  is convex. Since $\{S_t,\, t\in[0,T]\}$
is  a positive martingale, assuming enough integrablity,  $\{S_t^r,\, t\in[0,T]\}$ is a submartingale. Hence, for $t<t'$,
$$\E[S^r_t]\le \E[S^r_{t'}].$$
For $\epsilon>0$, such that   $u_+^*(t')-\varepsilon>1$, we have that
$$\E[S^{u_+^*(t')-\varepsilon}_t]\le \E[S^{u_+^*(t')-\varepsilon}_{t'}]<\infty,$$
and this implies $u_+^*(t)\ge u_+^*(t')-\varepsilon$, and the result follows from (\ref{sup-inf}).

For the negative abscissa, just observe that for $r<0$,  the same function  $\varphi(x)=x^{r}$  on $(0,\infty)$
is also convex, and apply the same reasoning. $\square$

\bigskip

\begin{figure}[htb]
\centering
%
%
%
\includegraphics[scale=0.6]{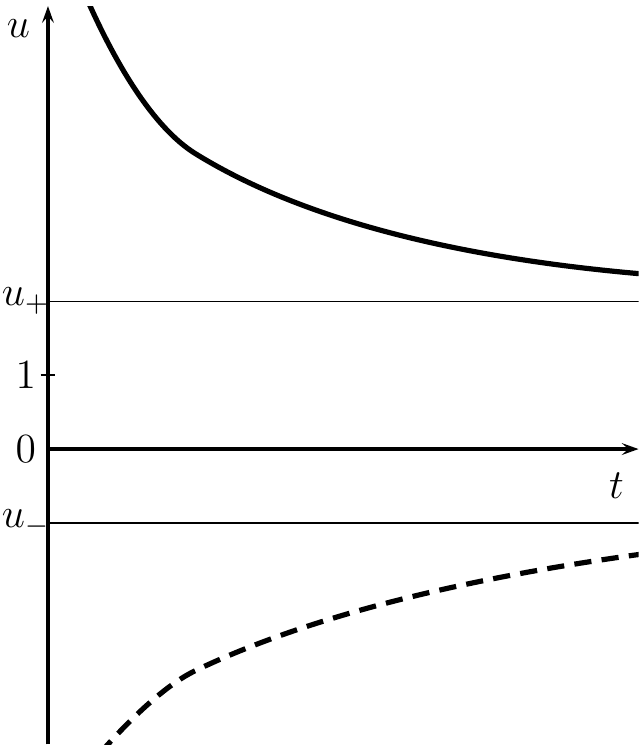}
\caption{Solid line: Plot of $u^*_+(t)$. Dashed line: Plot of $u^*_-(t)$.}
\label{pols-temps}
\end{figure}

For Heston model, we can be more precise.  From the bounds given in Theorem \ref{segonafita}, it is deduced
the behaviour of the abciss\ae\ for $t\to\infty$
$$\lim_{t\to\infty} u_+^*(t)=u_+ \quad\text{and}\quad \lim_{t\to\infty} u_-^*(t)=u_-,$$
and for $ t\to 0$:
$$\lim_{t\to 0} u_+^*(t)=\infty \quad\text{and}\quad \lim_{t\to 0} u_-^*(t)=-\infty.$$
A plot of these funcions is given is Figure \ref{pols-temps}.

\subsection{The effective vol-of-vol and the effective mean reversion}
 It is interesting that in
formula (\ref{umes}) the numbers $u_{\pm}$ only depend on the quotient $\omega:=a/c$.
Expressed in terms of this parameter
we have, for $\rho\ne \pm 1$,
$$u_{\pm}=\frac{1-2\omega\rho\pm\sqrt{4\omega^2+1-4\omega\rho}}{2(1-\rho^2)}.$$
Intuitively, this is a
consequence of the fact that the parameter $a$ (mean reversion) dampens the
stochastic volatility whereas $c$ (vol-of-vol) increases it, thus operating in
opposite directions.
However the parameters
$u_\pm$ correspond to the moments at infinite time, and in general the
relative strength of the parameters $a$ and $c$ in flattening the smile is time
dependent.  We propose to call this parameter the
{\it effective mean reversion factor},  and its inverse $c/a$ the
effective vol-of-vol.
The use of this parameter simplifies the study of $u_\pm$ for $\rho$ near $\pm 1$:
$$\lim_{\rho \to -1}u_-=-\frac{\omega^2}{2\omega+1}\quad \text{and} \quad \lim_{\rho \to -1}u_+=\infty,$$
and
\begin{equation*}
\lim_{\rho \to 1}u_-=\begin{cases}\dfrac{\omega^2}{2\omega-1},& \text{if\ $\omega\ < 1/2$}\\
\\
-\infty, & \text{if \ $\omega\ge 1/2$}
\end{cases}
\qquad \text{and} \qquad \lim_{\rho \to 1}u_+
=\begin{cases}\infty, & \text{if \ $\omega\le 1/2$}\\
\\
\dfrac{\omega^2}{2\omega-1},& \text{if \ $\omega\ > 1/2$}
\end{cases}
\end{equation*}

\section{Factorization of the moment generating  function of $\bs{X_t}$.}
\label{factorsec}

In this  section we will work exclusively with the random variable $X_t$ for $t>0$ fixed, and   the time $t$ will be
considered a parameter. Denote by ${\cal HL}(a,b,c,\rho,x_0,v_0,\mu,t)$ the law of $X_t$, that is,  a
probability on $\R$ that has the moment generating function given by (\ref{generating}). Observe that for all
$\lambda >0$,
$${\cal HL}(\lambda a,\lambda b,\lambda c,\rho,x_0,\lambda v_0,\mu,t/\lambda)={\cal HL}(a,b,c,\rho,x_0,v_0,\mu,t).$$
In particular,
 $${\cal HL}(a,b, c,\rho,x_0,v_0,\mu,t) =
{\cal HL}(at/2,bt/2, ct/2,\rho,x_0,v_0t/2,\mu,2).$$
Since we are interested in a property true for all parameters, it suffices to prove that property for arbitrary
$a,b,c,\rho,x_0,v_0,\mu,$ and $t=2$.
So, in all  proofs, we will take this value of $t$.

In this section we use some powerful theorems of complex variable analysis that we apply to the
complex  moment  generating function
\begin{align}
 \label{laplace2}
\Phi(z)=&
\exp\{x_0 z\}
 \ \bigg(\frac{e^{(a-c\rho z)t/2}}{\cosh(P(z)t/2)+(a -c\rho z)\sinh(P(z)t/2)/P(z)}\bigg)^{2a b/c^2}\notag \\
\notag \\
& \qquad \cdot\exp\bigg\{-v_0\ \frac{(z-z^2)\sinh(P(z)t/2)/P(z)}
{\cosh(P(z)t/2)+(a-c\rho z)\sinh(P(z)t/2)/P(z)}\bigg\},
\end{align}
defined in a neighborhood of zero, introduced in Subsection \ref{complexm}.
We will consider separately the different factors of this function.
and at a later stage, we will
combine them.

\subsection{The entire component}
Write
\begin{equation}
\label{fz}
F(z)=\cosh(P(z)t/2)+(a-c\rho z)\,\frac{\sinh(P(z)t/2)}{P(z)}.
\end{equation}
Recall that
$$p(z)=(a-\rho c z)^2+ c^2(z-z^2)$$
and $P(z)=\sqrt{p(z)}.$
As in Subsection \ref{extension}, but now in the complex plane, define the entire functions by the power series
\begin{equation}
\label{la}
L_1(z)=\sum_{n=0}^\infty \frac{z^n}{(2n)!},
\end{equation}
and
\begin{equation}
\label{lb}
L_2(z)=\sum_{n=0}^\infty \frac{z^n}{(2n+1)!}.
\end{equation}
We note that at each $z$,
$$L_1(z)=\cosh \sqrt{z} \quad \text{and}\quad L_2(z)=\frac{1}{\sqrt{z}}\, \sinh \sqrt z,$$
independently of the branch of the square root. Indeed, in every neighborhood that does not include zero, the previous
relations are true fixing  an arbitrary branch of the square root.   However,
 in the whole $\C$, the  functions $L_1$ and $L_2$ are defined by the power series and  not as a composition of $\cosh z$ or $\sinh z$ and
a particular branch of  $\sqrt z$. We consider the extension of $F(z)$ to an entire function using the functions
$L_1(z)$ and $L_2(z)$. Note that both $L_1(z)$ and $L_2(z)$ take real values on $\R$, and that $F(z)$
restricted to $\R$ coincides with the function $F(u)$ defined in (\ref{unificacio}).

The first interesting property of $F(z)$ is that it has all the zeroes real and simple. This can be deduced from a
 deep
theorem of Lucic \cite[Theorem 1]{Luc07} who  proves that all the singularities of the (complex) characteristic function of
 $X_t$ are purely imaginary.
In Appendix B there is  an alternative proof of this property. Specifically, we prove that

\begin{proposition} The zeroes of $F(z)$ are all real and simple.
\label{positives}

\end{proposition}

\bigskip

Using the Hadamard representation theorem (see the Appendix, Theorem B.\ref{hadamard-theo}), we deduce the following representation of
$F(z)$.

\begin{theorem}
\label{hadamard-fact}
 Let $\{a_n(t),\, n \ge 1\}$ be the   zeroes of $F(z)$. Then
\begin{equation}
\label{hadamard}
F(z)=e^{at/2}e^{\nu(t) z }\prod_{n=1}^\infty\Big(1-\frac{z}{a_n(t)}\Big)\exp\Big\{\frac{z}{a_n(t)}\Big\},
\end{equation}
where $\nu(t)\in \R$.
\end{theorem}

\begin{remark}

 The number $\nu(t)$ in (\ref{hadamard}) is determinated by
\begin{equation}
\label{constant-d}
F(1)=\cosh \big(\vert a-\rho c\vert t/2\big)+\text{sign}\,\{a-\rho c\}\, \sinh \big(\vert a-\rho c\vert t/2\big)
=e^{at/2}e^{\nu(t)  }\prod_{j=1}^\infty\Big(1-\frac{1}{a_j(t)}\Big)\exp\Big\{\frac{1}{a_j(t)}\Big\}.
\end{equation}

\end{remark}

\begin{remark}
In all of this Section  we exclude the case $\rho=1$ and $a=2c$ because the moment generating function simplifies and  has
only one pole in that case. Then the factorization is trivial.
\end{remark}

\subsection{The meromorphic component}
Now we will deal with the other factor of the function $\Phi(z)$ given in (\ref{laplace2}),
$$\frac{(z-z^2)\sinh(P(z)t/2)/P(z)}
{\cosh(P(z)t/2)+(a-c\rho z)\sinh(P(z)t/2)/P(z)}.$$
Write
$$G(z)=\frac{(1-z)\sinh(P(z)t/2)/P(z)}
{\cosh(P(z)t/2)+(a-c\rho z)\sinh(P(z)t/2)/P(z)},$$
where, as above, the function is extended to $\C$ using the entire functions $L_1(z)$ and $L_2(z)$. The function $G(z)$
is meromorphic  with simple poles at the roots of $F(z)$.
Using a theorem of Mittag-Leffler, see the Apendix, Theorem B.\ref{mittagtheo}, we prove

\begin{theorem}
\label{mittag-desc}
Let $\{a_n(t),\, n\ge 1\}$ be the zeroes of $F(z)$, ordered in the following way: $0<\vert a_1(t)\vert \le \vert a_2(t)
\vert \le \cdots$ .
  Then, for  $z\ne a_n(t)$, $\forall n\ge 1$,
\begin{equation}
\label{mittag}
G(z)=\frac{1}{2a}(1-e^{-at})-z\sum_{n=1}^\infty  \frac{b_n(t)}{a_n(t)^2(1-z/a_n(t))},
\end{equation}
where $b_n(t)>0$, and $\sum_n \frac{b_j(t)}{a_j^2(t)}<\infty$,
and the series converges uniformly in every compact set included in the disc $\vert z\vert <\vert a_1\vert$.

\end{theorem}

\subsection{Factorization of the  moment generating function of $\bs{X_t}$}
We return to the real moment generating function $M(u)$.
Combining the results of  the two preceding subsections we have that the  moment generating function $M(u)$ given
in (\ref{generating}) can be factorized as (we abreviate  $a_n(t)$ and $b_n(t)$ to $a_n$ and $b_n$ for a moment), for $u\in(-\vert a_1\vert,\, \vert a_1\vert),$
$$M(u)=\exp\{du\}\,
\prod_{n=1}^\infty\Big(1-\frac{u}{a_n}\Big)^{-2ab/c^2}
\exp\Big\{-\frac{2ab}{c^2}\frac{u}{a_n}+v_0 \frac{u^2 b_n}{a_n^2(1-u/a_n)}\Big\},
$$
where $d$ is a parameter to be defined  later. Each exponential can be written as
$$c_n u+ \frac{g_n u }{1-u/a_n}.$$
where
\begin{equation}
\label{parameters}
c_n=-(v_0 b_n+2ab/c^2)/a_n\quad\text{and}\quad g_n=v_0 b_n/a_n.
\end{equation}

\begin{remark}
When $4ab/c^2=k$ is a natural number,  we can identify the moment generating function as the one given by Janson
\cite[Theorem 6.2]{Jan97}, each eigenvalue with multiplicity   $k$. That means,
for such parameters, $X_t$ is in the (non-homogeneous)
second Wiener chaos. It is  well known that for such combination
of parameters the CIR model has the law of a sum of the squares of $k$ independent Ornstein-Uhlenbeck processes,
and this fact has been used for many applications. See, for example, Grasselli and Hurd \cite{GraHur05} and the references therein.

\end{remark}

On the other hand, for $u$ in a neighborhood of 0, the moment generating function of each factor
  can be written as
$$N_j(u)=e^{c_n u}\Big(1-\frac{ u}{a_n}\Big)^{-2ab/c^2}
\exp\Big\{ \frac{g_n u}{1-u/a_n}\Big\},$$
where $a_n\cdot g_n>0$.

\begin{proposition}
\label{inversio-simple} For $\xi>0$ and $\zeta,\gamma\in\R$, such that $\zeta\cdot\gamma>0$, the function
$$N(u)=\big(1-u/\gamma\big)^{-\xi}\exp\Big\{ \frac{\zeta u}{1-u/\gamma}\Big\}$$
for $u$ in a neighborhood of 0,
is the moment generating function  of an absolutely continuous random variable with probability density function
$$h(x)=\frac{1}{2}\bigg(\frac{x}{\zeta}\bigg)^{\tau}\exp{\{-(x+\zeta)\gamma\}}I_{2\tau}(2\gamma\sqrt{\zeta x})
\1_{(0,\infty)}(x),\ \text{if}\  \zeta,\gamma>0,$$
or
$$h(x)=\frac{1}{2}\bigg(\frac{x}{\zeta}\bigg)^{\tau}\exp{\{-(x+\zeta)\gamma\}}I_{2\tau}(-2\gamma\sqrt{\zeta x})
\1_{(-\infty,0)}(x),\  \text{if}\  \zeta,\gamma<0,$$
where $\tau=(\xi-1)/2$ and $I_{2\tau}$ is the Bessel function of index $2\tau$.

\end{proposition}

\noindent{\it Proof.}

For $\zeta,\gamma>0,$ the moment generating function $N(u)$ corresponds to the law of $Y_{1/(2\gamma)}$ for a  Bessel process
$$Y_t=\zeta+2\xi t+2\int_0^t\sqrt{Y_s}\, dW_s.$$
See Revuz-Yor \cite[Chapter 11]{RevYor99}.  Its density is  also  given
in Revuz-Yor \cite[page 441]{RevYor99}.

For $\zeta,\gamma<0,$ consider the random variable $Y$ defined above with parameters $-\zeta$ and $-\gamma$, and let
$Y'=-Y$. Its moment generating function is
$$M_{Y'}(u)=M_Y(-u)=\big(1-u/\gamma\big)^{-\xi}\exp\Big\{ \frac{\zeta u}{1-u/\gamma}\Big\}.$$
Its density is $h_Y(-x)$, where $h_Y$ is the density of $Y$. And the result follows. \quad $\square$

\bigskip

Return to the factors $N_j(u)$.  The term $e^{c_n u}$ is a translation factor of the random variable considered
above, so  $N_j(u)$ corresponds to a random variable with density $h_n(x-c_n),$ where $h_n$ is the probability density
function
given in Proposition \ref{inversio-simple}.

Finally, the product of characteristic functions of absolutely continuous laws with densities  $k_1$ and $k_2$
corresponds to the convolution
 $k_1 \star k_2$:
   $$k_1 \star k_2(x)=\int_{-\infty}^{+\infty} k_1(y)k_2(x-y)\, dy.$$
  Denote
by  $\bigstar_{j=1}^n k_j$ the product $k_1\star \cdots\star k_n$, that is  defined without ambiguity because
the convolution
 product is associative and commutative.

With all these ingredients we  construct a sequence of random variables that converges in law to $X_t$.
In general, the convergence in law does not imply the convergence of the corresponding probability density functions.
However, from a practical point of view, this fact does not matter.

\begin{theorem} The sequence of laws with densities $\bigstar_{j=1}^n \wt h_j$ converges to the law of $X_t+d(t)$,
where $\wt h_n(x)=h_n(x-c_n)$, and $h_n$ is the probability density function given in Proposition \ref{inversio-simple},
with parameters $$\xi=2ab/c^2, \,\gamma=a_n(t), \, \zeta=g_n(t),$$
where

\begin{enumerate}[1.]
\item $a_n(t)$ are the roots  of the function $F(u)$.
\item
$$b_n(t)=\frac{4p\big(a_n(t)\big)\big(1-a_n(t)\big)}{t p'\big(a_n(t)\big) p\big(a_n(t)\big)
-4c\rho p\big(a_n(t)\big)-p'\big(a_n(t)\big)(a-c\rho a_n(t))\big((a-c\rho a_n(t))t+2\big)}.$$
\item $c_n(t)$ and $g_n(t)$ are given in equation (\ref{parameters}).
\item $$d(t)=x_0-\frac{\rho ab t}{c}-\frac{2ab\nu(t)}{c^2}-\frac{v_0(1-e^{-at})}{2a}.$$
\item $\nu(t)$ is given in formula (\ref{constant-d}).
\end{enumerate}
\end{theorem}

\section*{Conclusions}

We have presented a new expression of the characteristic function of log-spot in Heston model that shows  its
good analytical properties and that facilitates its study. Through an analysis of the corresponding moment
generating function, we give numerical formulas   to obtain the absciss\ae\ of convergence, that have interesting
applications. As  examples we considered the computation of the parameters describing the asymptotic
behaviour of the volatility smile  for  extreme strikes, and the verification  that the model has enough moments
to price wing dependent deals.
 Another application may be the possibility of assessing the stability of the moments in
a  calibration of a Heston model.

In the second part of the paper, we factorized the moment generating function as an infinite product of Bessel type moment
generating function. This gives a new insight of  the Heston model, showing its complexity, and
its relationship with  Ornstein-Uhlenbeck processes.
Further,  each factor can be inverted, and a sequence of random variables that converges in law to log-spot
can be deduced.
 Though such sequence is not easy to manage, because it relies on the   computation of the  roots of a
(real) function and other parameters, and the  convolution of densities involving Bessel functions. The fact that
all computations are real can open the possibility of alternative methods to the numerical inversion of the characteristic
function that are currently used by practitioners.

\section*{Appendix}
\subsection*{Appendix A. Complex and real moment generating function}
For a sake of completeness, we recall some of the properties of the complex and real moment generating function
of an arbitrary random variable $X$.  We follow Hoffmann-J{\o}rgensen \cite{Hof94}. The complex function
$$M_\C(z)=\E[e^{zX}],$$
defined on the set
$${D}_\C=\big\{z\in \C: \ \E\big[\big\vert e^{zX}\big\vert\big]<\infty\big\},$$
is called the complex  moment generating function ({\itn cmgf} from now on) of $X$.
We will denote by $\inte_\C$ the interior of ${ D}_\C$.
The  restriction of $M_\C(z)$ to the real numbers is the moment generating function ({\itn mgf})
$$M_\R(u)=\E[e^{uX}],$$
defined on the set
$$D_\R=\big\{u\in \R: \ \E\big[ e^{uX}\big]<\infty\big\}.$$
(we also put $\inte_\R$ by the interior of $D_\R$).
For the present purposes it is convenient to maintain the double notation with the subindices $\R$ and $\C$.
Since for $z \in \C$ and $a\in \R$,
$ \vert \exp\{az\} \vert=\exp\{a\, \Re z\},$
we deduce that
$${D}_\C=\big\{z\in \C: \ \Re z \in {D}_\R\big\}=D_\R+i\,\R.$$
The most important property of both $M_\C(z)$ and $M_\R(u)$ is that they are analytic functions on the interior of its domains,
and that the Taylor expansion of $M_\C(z)$ in a point of the real axis has the same (real) coefficients than $M_\R(u)$.
Again, here, it is useful to introduce a double notation for the neighborhoods.
For $u\in \R$ and $r>0$, we denote by
 $B_\R(u,r)$ the neighborhood centered at $u$ with radius  $r$, and for $z\in \C$, $B_\C(z,r)$ is the
  neighborhood in  $\C$.
 If $D_\R\ne \{0\}$, the cmgf $M_\C$ (respectively the mgf $M_\R$)  is  analytic in  $\inte_\C$ (resp. $\inte_\R$).
and for  $z_0\in \inte_\C$ there is $r>0$ such that
\begin{equation}
\label{expansio}
M_\C(z)=\sum_{n=0}^\infty \,\frac{\E\big[ e^{z_0X}X^n\big]}{n!}\, (z-z_0)^n,\quad \forall z\in B_\C(z_0,r),
\end{equation}
and for  $u_0\in \inte_\R$ there is $r>0$ such that
\begin{equation}
\label{expansio2}
M_\R(u)=\sum_{n=0}^\infty \,\frac{\E\big[ e^{ u_0 X}X^n\big]}{n!}\, (u-u_0)^n,\quad \forall u\in B_\R(u_0,r).
\end{equation}
In particular,  if $0\in\inte_\R$, then $X$ has finite moments of all orders and, writing
$$m_n=E[X^n], \ n\ge 1,$$
we have that for some $r>0$,
$$M_\C(z)=\sum_{n=1}^\infty \frac{m_{n}}{n!}\, z^{n},\ \forall z\in B_\C(0,r),$$
and
$$M_\R(u)=\sum_{n=1}^\infty \frac{m_{n}}{n!}\, u^{n},\   \forall u\in B_\R(0,r).$$
Moreover,
$$\E[X^n]=M^{(n)}(0).$$

\begin{propositiona}
\label{coinc}
Assume that  there is a function  $\Phi(z)$ of the complex variable $z$ analytic in a neighborhood of \, 0 such that
$$M_\R(u)=\Phi(u),$$
for $u$ in a (real) neighborhood of 0. Then the function $\Phi(z)$ can be analytically continuated in
a strip that includes the imaginary axis, and   the characteristic function of $X$ is
$$\E[iuX_t]=\Phi(iu), \ \forall u\in \R.$$
\end{propositiona}

\noindent{\it Proof.} Since $M_\R(u)$ and $M_\C(z)$ have the same Taylor coefficients  in a neighborhood of 0, we have that
$\Phi(z)=M_\C(z)$ in a complex neighborhood of 0. Since $M_\C(z)$ is analytic in $\inte_\C$,
the function $\Phi$ can be continuated  to that strip that includes the imaginary axis.
 \quad $\square$

\bigskip

Remember that  $D_\R$ is an interval of $\R$ (it may be ${0}$), and the right extreme
(respectively the left extreme) of $D_\R$ is called the right--abscissa (resp. the left--abscissa)
of convergence.

\begin{propositiona}
\label{abscisa}
Assume that $0\in \inte_\R$, and that the right-abscissa of convergence, $\beta$, (respectively the left-abscissa $\alpha$)
is finite. Then $M_\R(u)$ has a singularity
at $\beta$ (resp. $\alpha$).
\end{propositiona}

\noindent{\it Proof.} When $X$ is positive, the cmgf coincides with the Laplace transform (except a change of sign on $z$)
of the distribution function of $X$. By  Widder \cite[Theorem II.5b]{Wid41}), $M_\C(z)$ has a singularity at the real point
$\beta$. From the fact that in the real axis both $M_\C$ and $M_\R$ have the same coefficients of the Taylor expansion, Widder's proof can be translated
to $M_\R$. For a general $X$, we can use the habitual technique of  decomposing a bilateral Laplace transform as the addition
of two unilateral ones (see Widder \cite[page 237]{Wid41})).
\quad $\square$

\subsection*{Appendix B. Proofs}

\noindent{\bf Proof of Proposition \ref{primerafita}}

We have that $p(0)=a^2$ and $p(1)=(a-\rho)^2\ge 0,$ then, given the form or $p(u)$,  it follows that  $1\le u_+$.
Hence, the  affirmation $[u_-,1]\subset D(X_t)$
 is obtained from the other points of the proposition.

\begin{figure}[htb]
\centering

%
%
%

\includegraphics[scale=0.8]{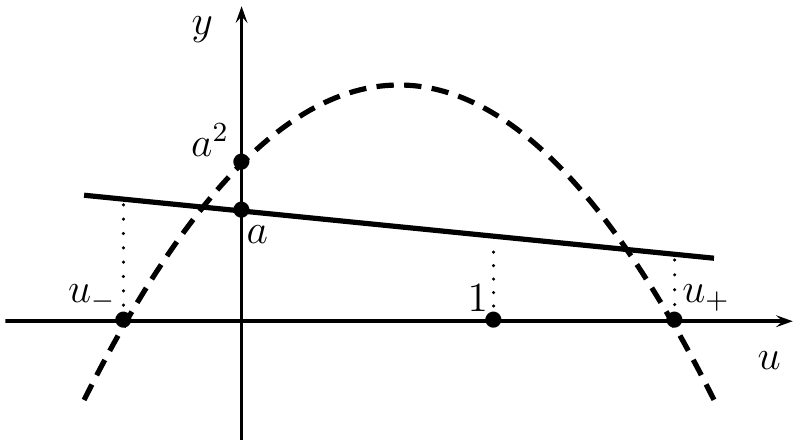}

\caption{Solid line: straightline $y=a-\rho c u$, case  $a>\rho c$, for $\rho>0$.
Dashed line: parabola $p(u)$}
\label{recta-parabola}
\end{figure}

{\bf 1.} Consider   the case $a\ge \rho c$, $\rho\ne \pm 1.$
A bit of algebra shows that $a-\rho c u_->0$ and $a-c\rho u_+>0$. So the straightline
$y=a-\rho c u$ is positive for $u\in[u_-,u_+]$ (see Figure \ref{recta-parabola}). Since $\cosh x\ge 1, \forall x\in \R$, and $\sinh x>0, \forall x>0$,
it follows that $f(u)>0$ in $[u_-,u_+]$. Note that if $a=\rho c$, then $u_+=1$.

\medskip

When $\rho=-1$, then  $u_-=-a^2/(c(c+2a))<0,$ and it is trivial that   $a+c u_->0$. So  the straithgtline $y=a+c u$ is positive for $u\ge u_-$, and hence
$f(u)>0$ in $[u_-,\infty)$.

When $\rho=1$, and $c\ne 2a$,  then $u_+=a^2/(c(2a-c))$, and some calculations shows that
in this case ($a\ge c)$, then
$$\frac{a}{c}\ge \frac{a^2}{c(2a-c)},$$
and hence, $a-cu_+\ge 0$. Then   $a-cu\ge 0$ for all $u\le u_+$, and it follows $f(u)>0$ in $(-\infty,u_+].$

\medskip

{\bf 2.} Consider the case $a <\rho c$ (note that this implies $\rho >0$).
Then $a-c\rho u_->0$ and $a-c\rho \big(\frac{a}{c\rho}\big)=0$. So the straighline $y=a-c\rho u$ is positive for
$u\in[u_-,a/(c\rho)]$. Thus $f(u)$ has no zeroes in such interval. On the other hand,  the real zeroes of $f(u)$ in $[u_-,u_+]$
 are the same of the ones corresponding
to the  function
$$1+\frac{a-c\rho u}{P(u)}\tanh(P(u)t/2),$$
because of the real hyperbolic cosine is never zero.
For $u\ge  a/(\rho c),$ and the bound $0<\tanh x < 1$, for all $x>0$, we have
$$1+\frac{a-c\rho u}{P(u)}\tanh(P(u)t/2)>1+\frac{a-c\rho u}{P(u)},$$
and when  $u\in [a/(\rho c),1]$, from $c^2u(1-u)\ge 0$ it is deduced that
$$1+\frac{a-c\rho u}{P(u)}\ge 0,$$
hence  $f(u)>0$ also in $[a/(\rho c),1]$. Now,  we  analyze the behaviour of $f(u)$ for $u\in[1, u_+]$.
Assume $\rho\ne 1$.
From $$\lim_{u \nearrow u_+ }\frac{\sinh\big (P(u)t/2\big)}{P(u)t/2}=1.$$
it follows that
$$\lim_{u\nearrow   u_+ }   f(u)=1+(a-c\rho u_+)t/2.$$
Consider the straightline $y=a-c\rho u$. We have $y(0)=a$ and $y(1)=a-c\rho<0$, and since $u_+>1$,
 then $a-c\rho u_+<0$. Let $t_0=2/(c\rho u_+-a)>0.$

\begin{enumerate}[(i)]
\item If $t< t_0$,  from $(c\rho u_+-a)t_0/2=1$, then for all $u\in[1,u_+],$
$$0<(c\rho u-a)t_0/2<1, $$
and thus
$$0<(c\rho u-a)t/2<1.$$
Hence, $\forall u\in[1,u_+],$
$$\frac{(c\rho u-a)t}{2}\ \frac{\sinh \big(P(u)t/2\big)}{P(u)t/2}<\frac{\sinh \big(P(u)t/2\big)}{P(u)t/2}\le \cosh \big(P(u)t/2\big),$$
because  ${\sinh x}/{x}\le \cosh x $. Then $ f(u)>0.$

\item Let $t\ge t_0$. Then $\ f(1)>0$ and $ f(u_+)<0,$ and thus there is at least one zero of
 $f(u)$ in $[1,u_+].$ The fact that there is only one zero is proved in Section \ref{factorsec}.
 \end{enumerate}

Finally, the case $\rho=1$ and $a<c$\,  or $a=2c$ are studied in a similar way. $\square$

\bigskip

\bigskip

\noindent{\bf \large Proof of  the results of Section   \ref{factorsec}}.

\medskip

For easy reference, we recall here the two main theorems used in the proofs, expressed in the form that we need.
The first one is the factorization Theorem of Hadamard (see, for example,  Titchmarsch \cite{Tit52})

\begin{theoremb}
\label{hadamard-theo}
Let $F(z)$ an entire function, $F(0)\ne 0$,  with  roots $a_1, \,a_2,\dots,$  such that
$$\sum_{n}1/\vert a_n\vert^2<\infty.$$
Then $F$ can be represented as
$$F(z)=F(0)e^{C z}\prod_{n} \Big(1-\frac{z}{a_n}\Big)e^{z/a_n}.$$
\end{theoremb}

The second theorem is due to Mittag-Leffler (Titchmarsch \cite[page 110]{Tit52})
\begin{theoremb}
\label{mittagtheo}
Let $G(z)$ be a meromorphic function such all the poles are simple, denoted by  $a_1, \,a_2,\dots,$ where
$0<\vert a_1\vert \le \vert a_2\vert \le \cdots,$ and with residues
at the poles  $b_1,\,b_2\dots$ respectively. Suppose that there is a sequence of closed contours $C_n$ such that
$C_n$ includes $a_1,\dots, a_n$ but no other poles, such
that the minimum distance $R_n$ of $C_n$ to the origin tends to infinite with n, while the length of $C_n$ is $O(R_n)$,
and on $C_n$, $f(z)=o(R_n)$.  Then
$$G(z)=G(0)-z\sum_{n}\frac{b_n/a_n^2}{1-z/a_n}.$$
\end{theoremb}

In order to prove Proposition  \ref{positives} we need three lemmas.

\begin{lemmab}
\label{fita1}
 For all $z\in \C$,
$$\big \vert \coth \sqrt z \big\vert\le \sup_{\vert w\vert= \sqrt{\vert z\vert}} \big \vert \coth w \big\vert,$$
and similarly for $\tanh \sqrt z$.
\end{lemmab}

\medskip

\noindent{\it Proof.}

Consider  the formula
\begin{equation}
\label{fitcot}
\vert \coth z \vert^2=\frac{\sinh ^2x+\cos^2y}{\sinh ^2x+\sin^2y}.
\end{equation}
 Write $z=r\,e^{i\theta}$ and choose an arbitrary branch of the square root, for
example, take the principal one. Then
\begin{align*}
\big \vert \coth \sqrt  z \big\vert ^2&=\frac{\sinh^2(\sqrt r\, \cos(\theta/2))+\cos^2(\sqrt r\, \sin(\theta/2))}
{\sinh^2(\sqrt r\, \cos(\theta/2))+\sin^2(\sqrt r\, \sin(\theta/2))} \\
&\le \sup_{\phi\in[0,\,2\pi)}
\frac{\sinh^2(\sqrt r\, \cos \phi)+\cos^2(\sqrt r\, \sin \phi)}
{\sinh^2(\sqrt r\, \cos \phi)+\sin^2(\sqrt r\, \sin \phi)} \\
&= \sup_{\vert w\vert =\sqrt r}\big \vert \coth   w \big\vert ^2. \quad \square
\end{align*}

\begin{lemmab}
\label{fita2}
For every $\varepsilon>0$ there is  $n_0\ge 1$ such that for $n\ge n_0$,
$$\sup_{\vert z\vert=(n+\frac{1}{2})\pi}\vert \coth z\vert <1+\varepsilon$$
and
$$\sup_{\vert z\vert=n\pi}\vert \tanh z\vert <1+\varepsilon.$$
\end{lemmab}

\smallskip

\noindent{\it Proof}

\begin{figure}[hbt]
\centering
%
%

\includegraphics[scale=0.4]{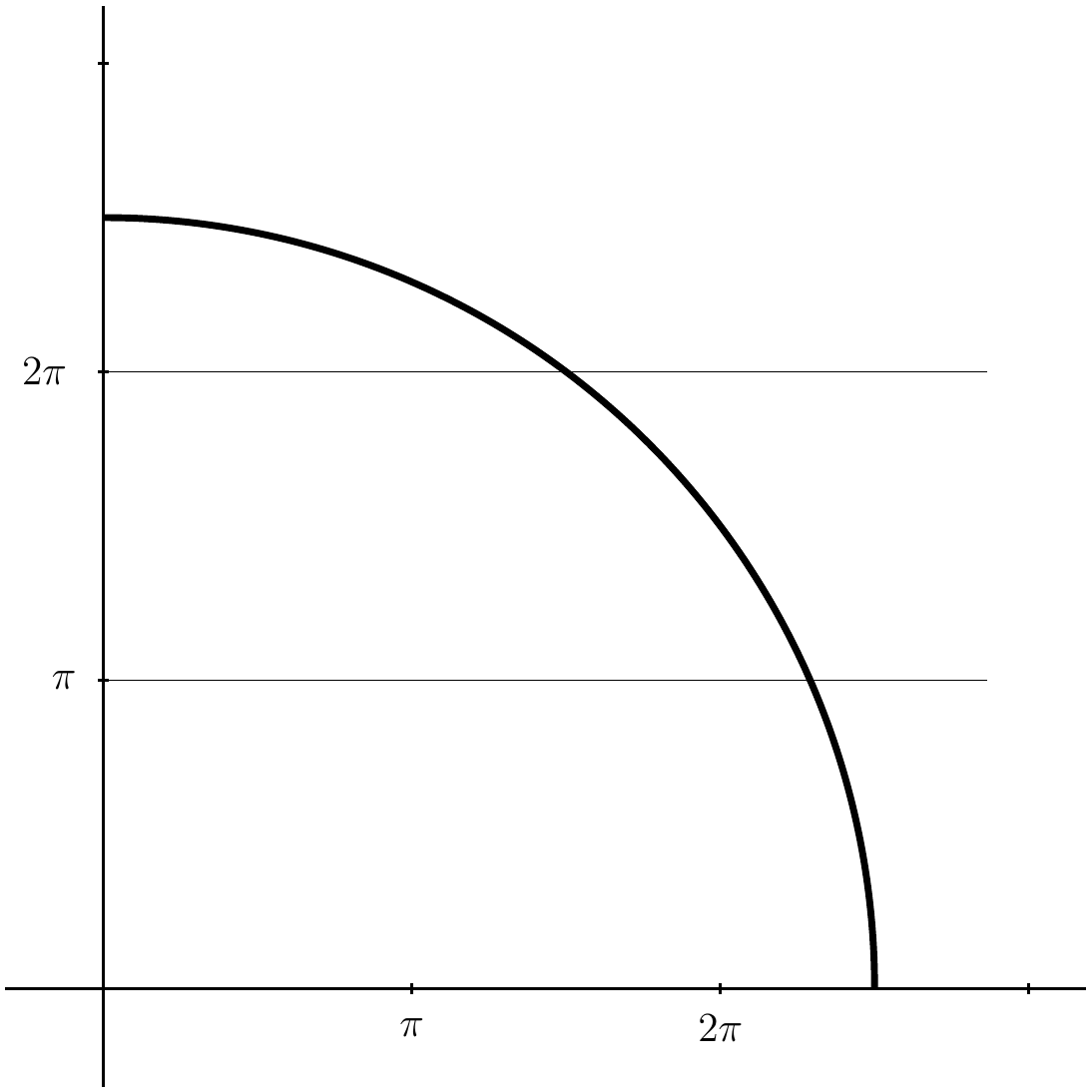}


\caption{Arc  $\vert z\vert=(n+\frac{1}{2})\pi, \ \Re z\ge 0,\ \Im z\ge 0$ for $n=2$.}

\label{circum}
\end{figure}

From the formula (\ref{fitcot})
it is clear that  $\vert \coth z\vert$ is the same for the  points    $z= x+iy, \, -x+iy,\ x-iy,\ -x-iy$.
Hence, to bound $\coth z$ in the circle
  $\{z\in \C: \ \vert z\vert=(n+\frac{1}{2})\pi\}$ we can restrict ourselves  to study the arc with
 $\Re z \ge 0$ i $\Im z\ge 0.$ (see the
 Figure \ref{circum} for $n=2$).

\begin{figure}[htb]
\centering
\subfigure[]{%
\includegraphics[width=6cm]{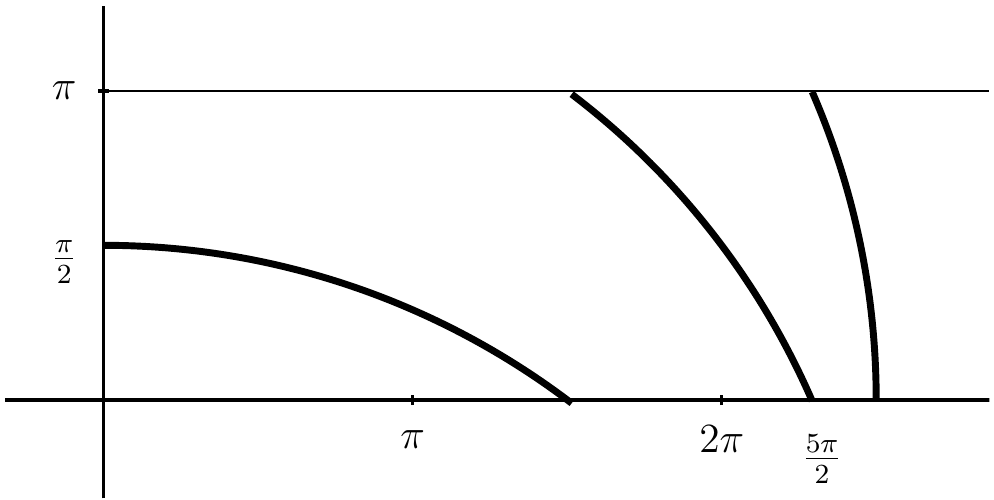}}
\hspace{2cm}
\subfigure[]{%
\includegraphics[width=6cm]{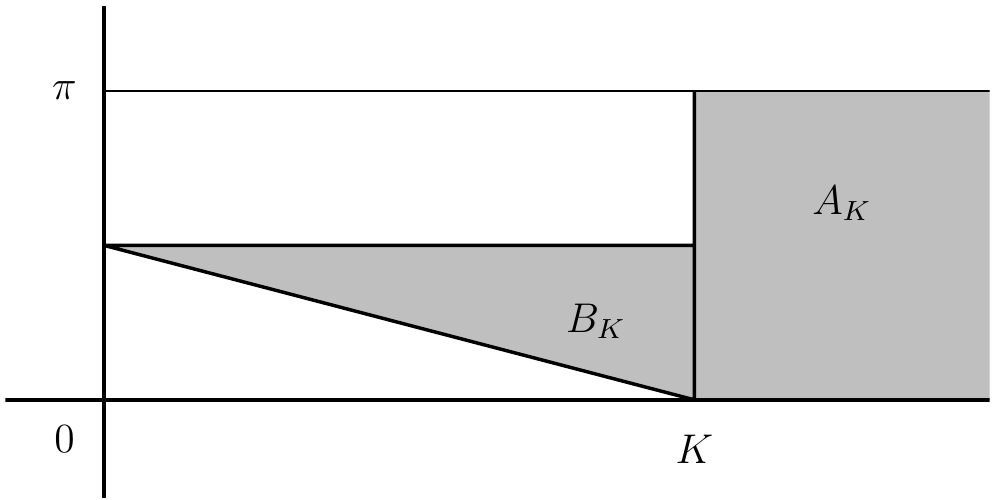}}
\caption{(a) Translation to the arc of Figure \ref{circum} to the strip $0\le \vert \Im z\vert\le \pi$.
(b) The translation for all $n\ge n_0$ are included in the shadow region}
\label{regio}
\end{figure}

Given the periodicity  $\coth(z+k\pi i)=\coth z$, for all  $k\in \Z$, it suffices to bound the translation
to the strip
  $0\le \Im z \le \pi$  (see Figure \ref{regio} (a)).
For $K$ big enough, for  $n\ge n_0$ all translations are included in the shaded region
$A_K\cup B_K$
of Figure \ref{regio} (b). Thus, it suffices to prove that for $K$ big enough,
$\sup_{z\in A_K\cup B_K}\vert \coth z \vert <1+\varepsilon.$
 For $z=x+iy\in A_K$, (that is, $ x  >K$), we have
$$\vert \coth z \vert^2=\frac{\sinh ^2x+\cos^2y}{\sinh ^2x+\sin^2y}
\le \frac{\sinh ^2x+1}{\sinh ^2x}=1+\frac{1}{\sinh^2x}\le 1+\frac{1}{\sinh^2 K},$$
that goes to 1 when   $K\to\infty.$

Now we study the bound on the triangle  $B_K$ of the Figure \ref{regio} (b).
By the principle of the maximum, we need only  to study the function on the border of $B_K$. On the vertical side,
$z=K+iy, \ y\in[0,\pi/2]$, works the bound that we have found above.
On the horizontal side,  $z=x+i\frac{\pi}{2},\ x\in[0,K]$, by the formula (\ref{fitcot}),
$$\vert \coth z \vert^2=\frac{\sinh ^2x}{\sinh ^2x+1}\le 1.$$
Finally, the hypotenuse is
 $z=x+iy$ with $y\in [0,\pi/2]$ i
$x=-2K y/{\pi}+K.$
For  $y\in[\pi/4,\, \pi/2],$ we have that   $\cos 2 y\le 0$, and by
$$\vert \coth z \vert^2=\frac{\cosh 2x+\cos 2y}{\cosh 2x-\cos 2y},$$
we deduce that
$$\vert \coth z \vert\le 1.$$
For  $y\in [0,\, \pi/4]$, we have
$x\in[K/2,K]$, and again by  (\ref{fitcot}) we obtain the bound.

\bigskip

To bound  $\vert \tanh z\vert $ for  $\vert z\vert=n\pi$ we do the same reductions as before, and use the relationship
$$\tanh z=\coth\big(z-i\frac{\pi}{2}\big),$$
to translate the bounds of $\coth z$ to $\tanh z$. \quad $\square$

\bigskip

\begin{figure}[htb]
\centering
\includegraphics[width=6cm]{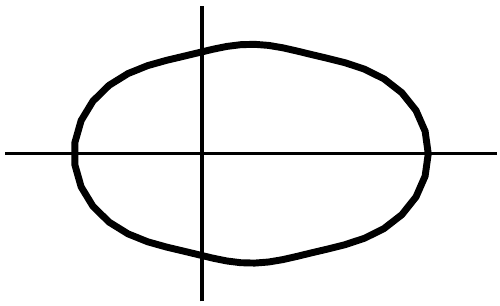}
\caption{Contour $\vert p(z)\vert =d$, for $d$ big.}
\label{contorn}
\end{figure}

 We also  need the properties of the contour determinated by  the polynomial $p(z)$.
For    $d>0$ let  $C_d$ be the contour
$$C_d=\{z\in\C:\, \vert p(z)\vert =d\},$$ and denote by   $L_d$ its length   and  by  $R_d$ the minimum distance from
$L_d$ to the origin. See Figure \ref{contorn}

\begin{lemmab}
\label{lemma-mittag}
With the above notations, for $d$ big enough,  $C_d$ is a homotopic to 0, $\lim_{d\to \infty} R_d=\infty,$
and  $L_d \sim o( R_d^2)$, when  $d\to\infty$.
\end{lemmab}

\smallskip

\noindent{\it Proof}

Since  $p(z)$ is a  second degree polynomial, there are two constants $K_1,\, K_2>0$ such that
for $\vert z \vert$ big enough,
$$K_1\vert z\vert ^2\le \vert p(z)\vert \le K_2\vert z\vert ^2.$$
It follows that for $d$ big enough, the circles with radius
 $d/K_2$ i $d/K_1$ bound lower and upper the contour $C_d$.
 This implies that $C_d$ is closed.

 The polynomial $p(z)$ has real coefficients and has a positive and a negative root ($\rho\ne \pm 1$), and
 since the length of a contour does not change by translation, we can consider that
 $p(z)=z^2-a$, with $a>0$ (change also $d$ in order that the coefficient of $z^2$ is 1). Write $z=x+iy$; then the equation that determines $C_n$ is
 $$(x^2+y^2+a)^2-4ax^2=d^2.$$
 In polar coordinates, $x=r(\theta)\cos \theta$ and $x=r(\theta)\sin \theta$, the contour  is given by
 $$r^4(\theta)-2 a \cos(2\theta) r^2(\theta)+a^2-d^2, \ \theta\in[0,2\pi].$$
 For $d$ big enough, we need only to consider the solution
 $$r(\theta)=\sqrt{a\cos(2\theta)+\sqrt{a^2cos(4\theta)+d^2}},$$
 that determines a curve clearly homotopic to 0. To compute the length $L_d$, by symmetry,
 $$L_d=4\int_0^{\pi/2}\sqrt{ r^2(\theta)+\big(r'(\theta)\big)^2}\, d\theta.$$
 We have $\lim_{d\to\infty} r(\theta)/d=0$ and $\lim_{d\to \infty} r'(\theta)/d=0$, uniformly in $\theta\in[0,\pi/2]$,
and  the lemma follows. \quad $\square$

\bigskip

\noindent{\bf Proof of Proposition \ref{positives}}

Remember that we consider $t=2$. First, note that for $\rho\ne\pm 1$,
 $$\lim_{z\to\infty}\frac{(a-c\rho z)^2}{(a-c\rho z)^2+c^2(z-z^2)}=\frac{\rho^2}{\rho^2-1}.$$
Hence,
\begin{equation}
\label{limit}
\lim_{z\to\infty}\Bigg\vert \frac{ a-c\rho z}{P(z)}\Bigg\vert
=\frac{\vert \rho\vert}{\sqrt{1-\rho^2}}.
\end{equation}
Consider the following  three cases:
\smallskip

\noindent{\bf Case 1.}  $\vert \rho\vert/{\sqrt{1-\rho^2}}>1.$

\medskip

\noindent{\bf First step.} The objective of this step is to prove that in the contour
$$ C_n=\{z\in C:\, \vert p(z)\vert=(n+0.5)^2\pi^2\},$$  for $n$ big,
\begin{equation}
\label{desigualtat1}
\vert\cosh P(z)\vert < \Big\vert (a -c\rho z)\frac{\sinh P(z)}{P(z)}\Big\vert.
\end{equation}
Then, by Rouche Theorem, both functions  $F(z)$ and  $(a -c\rho z)\sinh(P(z))/P(z)$ have the same number of
zeroes in $C_n$. To prove the inequality (\ref{desigualtat1}),  take $\varepsilon>0$
 such that $1<1+ \varepsilon <\vert \rho\vert/{\sqrt{1-\rho^2}}.$
By Lemmas {B3} and {B4},
 for $n$ big enough, and $z\in C_n$, we have
$\vert \coth P(z)\vert < 1+\varepsilon,$ and apply the limit (\ref{limit}).

\medskip

\noindent{\bf Second step.} Study of the zeroes of $F(z)$. Assume  that $\rho\ne 0$.
The zeroes of   $(a -c\rho z)\sinh(P(z))/P(z)$ are
 $z=a/(c\rho)$ and the zeroes of $\sinh(P(z))/P(z)$. We have that
$$\frac{\sinh P(z)}{P(z)}=\prod_{n=1}^\infty\Bigg(1+\frac{p(z)}{n^2\pi^2}\Bigg).$$
Hence, the zeroes of this function are the roots of the second order polynomial
 $$p(z)=-n^2\pi^2, \ n=1,\, 2,\dots,$$
 that we denote by $\alpha_{\pm n}$;
 write also   $u_-=\alpha_{-0}$ and $u_+=\alpha_{+0}$, in agreement with the notations of Section
 \ref{domini}.
They all  are  real and
 $$\cdots\alpha_{-(n+1)}<\alpha_{-n}<\cdots\alpha_{-1}<\alpha_{-0}<0<\alpha_{+0}<\alpha_{+1}<\cdots<\alpha_{+n}<\alpha_{+(n+1)}<\cdots.$$
So, for  $n$ big, $F(z)$ has $2n+1$ zeroes in $C_n$.

Now we  count the number of real roots of $F(z)$ in $C_n$. Assume that $a>c\rho$ and $\rho>0$. We saw in the proof of    Proposition \ref{primerafita}  that $a/(c\rho)>u_+$, and
assume that $a/(c\rho) \in (\alpha_{+k}, \alpha_{+(k+1)})$ (the case  $a/(c\rho)=\alpha_{+k}$  needs to be studied as a
particular case). Then,
\begin{enumerate}[(i)]
\item  In each interval $(\alpha_{+n}, \alpha_{+(n+1)})$, for $n=0,\dots k-1$,  $F(u)=\wt f(u)$, where
$\wt f(u)$  was defined in (\ref{funcioq}), and  $\wt f(u)$ has one root. This is deduced because the roots of
$\wt f(u)$ in such intervals are the solutions of
\begin{equation*}
 \tan \wt P(u)=-\frac{\wt P(u)}{a-c\rho u}.
 \end{equation*}

\item  In $(\alpha_{+k}, \alpha_{+(k+1)})$, $\wt f(u)$ has  $2$ roots. This claim is proved, observing  that
$a-c\rho \,\alpha_{+k}>0$ and $a-c\rho \,\alpha_{+k}<0$, and the the curve $-\wt P(u)/(a-c\rho u)$ cuts two times
the curve  $\tan \wt P(u)$ in that interval.

\item In each $(\alpha_{-(n+1)},\alpha_{-n}),\ n\ge 0$, the function $\wt f(u)$ has one root. This is proved as in point (i).

\end{enumerate}
All the other possibilities for $a,\, c$ and $\rho$ are discussed in a similar way, and we obtain that $F(z)$ has at least $2n+1$ real roots in $C_n$.
So the Theorem follows.

\bigskip

\noindent{\bf Case 2.} $\vert \rho\vert/{\sqrt{1-\rho^2}}<1.$ Here use the contours
$$ C_n'=\{z\in C:\, \vert p(z)\vert=n^2\pi^2\},$$
and prove that for $z\in C'_n, $ ($n$ big),
\begin{equation}
\label{desigualtat}
  \Big\vert (a -c\rho z)\frac{\sinh P(z)}{P(z)}\Big\vert<\vert\cosh P(z)\vert.
\end{equation}
and finish the proof as in Case 1.

\bigskip

\noindent{\bf Case 3.} $\rho=\pm\sqrt 2/2$.
 Write $\rho_n=\rho+1/n$, for $n\ge 4$, and let $F_n$ be the function $F$ with $\rho$ changed by $\rho_n$.
 We have  $F_k(z)\to F(z)$ as $k\to \infty$, uniformly in every disc. By Hurwitz theorem (see Titchmarsch \cite{Tit52}), the roots of $F(z)$ in such disc
 are the limit points of the roots of $F_k(z)$ in the disc. So the roots of $F(z)$ are also real. \quad $\square$.

\bigskip

\noindent{\bf Proof of Theorem \ref{hadamard-fact}}

In the proof of Proposition \ref{positives} we saw that for $n$ big, positive or negative, in each interval
$(\alpha_{+n}, \alpha_{+(n+1)})$ and $(\alpha_{-(n+1)}, \alpha_{-n})$ there is one and only one    root of $F(z)$, where
$$p(\alpha_{\pm n})=-n^2\pi^2.$$
Since $p(u)$ is a second degree polynomial, there are constants $K_1,\, K_2>0$
such that for $\vert u \vert$ big enough,
$$K_1\,  u ^2\le \vert p(u)\vert \le K_2\, u ^2,$$
and then
\begin{equation}
\label{fita-arrels}
\frac{1}{\sqrt {K_2}}n \pi \le \vert \alpha_{\pm n}\vert\le \frac{1}{\sqrt {K_1}}n \pi.
\end{equation}
Hence,
$$\sum_n \frac{1}{a_n^2}<\infty \quad \text{and}\quad \sum_n \frac{1}{\vert a_n\vert }=\infty.$$
By Hadamard factorization  Theorem B.\ref{hadamard-theo}
  we obtain the
representation  (\ref{hadamard}).
$\square$

\begin{remarkb}
An alternative way to prove the previous result is using that the order (as entire function) of $F(z)$ is less or equal
to 1. This can be deduced from that fact that
by the definitions (\ref{la}) and (\ref{lb})
$$\vert L_j(z)\vert\le \cosh (\sqrt{\vert z\vert}), \ j=1,2,$$
and that the order of $L_1(z)$ can be easily computed from  (\ref{la}) and   is 1/2.
\end{remarkb}

\bigskip

\noindent{\bf Proof of Theorem \ref{mittag-desc}}

In a similar way that  in the proof of Proposition \ref{positives},  we are going to prove that there is a constant
$C>0$ such that  for $z\in C_n$ or $z\in  C_n'$, for $n$ big enough,
$$\big\vert G(z) \big\vert <C,$$
where $C_n$ and $C'_n$ where defined in Theorem \ref{positives}. Then, thanks to  Lemma B.\ref{lemma-mittag}, we can apply
 the
theorem of Mittag-Leffler
B. \ref{mittagtheo}
 that gives the expression (\ref{mittag}).  Consider three cases:

\medskip

\noindent{\bf Case 1.} $\vert \rho\vert/{\sqrt{1-\rho^2}}>1.$ Take $\varepsilon>0$ such that $1<1+ \varepsilon <\vert \rho\vert/{\sqrt{1-\rho^2}}.$
By Lemmas B3 and B4,
for $n$ big enough, and $z\in C_n$,
$\vert \coth P(z)\vert < 1+\varepsilon.$ Therefore
\begin{align*}
\big\vert G(z)\big\vert& =
\Bigg\vert \frac{(1-z)/P(z)}
{\coth P(z)+(a-c\rho z)/P(z)}\Bigg\vert
 \le
\frac{\big\vert(1-z)/P(z)\big\vert}{\Big\vert\big\vert\coth P(z)\big\vert-\big\vert(a-c\rho z)/P(z)\big\vert\Big\vert}\\
&= \frac{\big\vert(1-z)/P(z)\big\vert}{\big\vert(a-c\rho z)/P(z)\big\vert-\big\vert\coth P(z)\big\vert}
<C.
\end{align*}

\bigskip

\noindent{\bf Case 2.}  $\vert \rho\vert/{\sqrt{1-\rho^2}}<1.$
Let $\delta>0$ be such that $\vert \rho\vert/{\sqrt{1-\rho^2}}<\delta<1,$
and $\varepsilon>0$ such that $(1+\varepsilon)\delta<1$. Then, for $z\in \wt C_n$ ($n$ large),
\begin{align*}
\big\vert G(z)\big\vert  =
\Bigg\vert \frac{(1-z)\tanh P(z)/P(z)}
{1+(a-c\rho z)\tanh P(z)/P(z)}\Bigg\vert
\le \frac{\big\vert \tanh P(z) \big\vert\, \big\vert(1-z)/P(z)\big\vert}
{ 1-\big\vert\tanh P(z)\big\vert\,\big\vert(a-c\rho z)/P(z)\big\vert}
<C.
\end{align*}

\bigskip

In both cases 1 and 2, by Mittag--Leffler Theorem
B. \ref{mittagtheo},
\begin{equation}
\label{mittag-2}
G(z)=G(0)-z\sum_{n=1}^\infty  \frac{b_n}{a_n^2(1-z/a_n)},
\end{equation}
where $b_n$ is the residue of $G(z)$ in the pole $a_n$. Since $G(z)$ is the quotient of two entire functions,
and the pole is simple,
$$b_n=\frac{(1-a_n)\sinh(P(a_n))/P(a_n)}{F'(a_n)}.$$
using that $a_j$ is a root of $F(z)$,
 differentiating and simplifying we obtain that
the corresponding residue is
$$b_n=\frac{2p(a_n)(1-a_n)}{p'(a_n) p(a_n)-2c\rho p(a_n)-p'(a_n)(a-c\rho a_n)(a-c\rho a_n+1)}.$$
Consider the function defined by the previous relation:
\begin{equation}
\label{funcio-pol}
g(u)=\frac{2p(u)(1-u)}{p'(u) p(u)-2c\rho p(u)-p'(u)(a-c\rho u)(a-c\rho u+1)}.
\end{equation}
Then
\begin{equation}
\label{limit-pol}
\lim_{u\to\pm \infty}g(u)=
\begin{cases} \frac{1}{c^2}, & \text{if $\rho^2\ne 1$}\\
 \frac{2}{c^2}, & \text{if $\rho^2= 1$ and $c\ne 2 a \rho$.}
 \end{cases}
 \end{equation}
So, for large $n$, it is clear that $b_n>0$. For small $n$ the positivity of $b_n$ is proved through an analysis of the
 sign of $f(u)$ and $f'(u)$ in the different intervals where there are located the  roots of $f$.
 These roots and its location are clearly studied in  Lucic \cite{Luc07}.

\bigskip

\noindent{\bf Case 3.}  $\rho=\pm\sqrt 2/2$. As in case 3 of  Theorem \ref{hadamard-fact}, the result  is obtained by continuity,
using monotone convergence Theorem. $\square$

\bigskip

\noindent{\bf\large Acknowledgements.} We would like to acknowledge professor Daniel Dufresne, from Melbourne University,
 from whom we learned the ideas about the moment generating function expressed in Lemma \ref{extensio}. We also are very grateful
 to professors  Armengol Gasull and Joan Josep Carmona from the Maths Department of the Universitat Aut{\`o}noma de Barcelona
 for helpful conversations.

\bibliographystyle{abbrv}
\bibliography{bibliografia}

\begin{thebibliography}{10}

\bibitem{AlbHanMaySch07}
H.~Albrecher, P.~Mayer, W.~Schoutens, and J.~Tistaert.
\newblock The little {H}eston trap.
\newblock {\em Wilmott Magazine}, pages 83--92, January, 2007.

\bibitem{AndPie07}
L.~B.~G. Andersen and V.~V. Piterbarg.
\newblock Moment explosions in stochastic volatility models.
\newblock {\em Finance Stoch.}, 11(1):29--50, 2007.

\bibitem{CoxIngRos85}
J.~C. Cox, J.~E. Ingersoll, Jr., and S.~A. Ross.
\newblock A theory of the term structure of interest rates.
\newblock {\em Econometrica}, 53(2):385--407, 1985.

\bibitem{Ban08}
S.~del Ba{\~n}o~Rollin.
\newblock Spot inversion in the {H}eston model.
\newblock {\em CRM Prepint 837}, 2008.

\bibitem{Duf01}
D.~Dufresne.
\newblock The integrated square-root process.
\newblock {\em Research Collections (UMER). Preprint}, 2001.
\newblock http://repository.unimelb.edu.au/10187/1413.

\bibitem{Duf08}
D.~Dufresne.
\newblock The distribution of realized volatility in stochastic volatility
  models.
\newblock {\em Preprint}, 2008.

\bibitem{Fel51}
W.~Feller.
\newblock Two singular diffusion problems.
\newblock {\em Ann. of Math. (2)}, 54:173--182, 1951.

\bibitem{Gat04}
J.~Gatheral.
\newblock A parsimonious arbitrage-free implied volatility parametrization with
  application to the valuation of volatility derivatives.
\newblock {\em Presentation at Global Derivatives \& Risk Management Madrid
  2004}, 2004.
\newblock www.math.nyu.edu/fellows fin math/gatheral/madrid2004.pdf.

\bibitem{Gat05}
J.~Gatheral.
\newblock {\em The volatility surface. {A} practitioner's guide}.
\newblock Wiley, New York, 2005.

\bibitem{GraHur05}
M.~R. Grasselli and T.~R. Hurd.
\newblock Wiener chaos and the {C}ox-{I}ngersoll-{R}oss model.
\newblock {\em Proc. R. Soc. Lond. Ser. A Math. Phys. Eng. Sci.},
  461(2054):459--479, 2005.

\bibitem{Hes93}
S.~Heston.
\newblock A closed form solution for options with stochastic volatility with
  applications to bond and courrency options.
\newblock {\em Review of Financial Studies}, pages 327--343, 6, 1993.

\bibitem{Hof94}
J.~Hoffmann-J{\o}rgensen.
\newblock {\em Probability with a view toward statistics. {V}ol. {I}}.
\newblock Chapman \& Hall Probability Series. Chapman \& Hall, New York, 1994.

\bibitem{Jan97}
S.~Janson.
\newblock {\em Gaussian {H}ilbert spaces}, volume 129 of {\em Cambridge Tracts
  in Mathematics}.
\newblock Cambridge University Press, Cambridge, 1997.

\bibitem{Lee04}
R.~W. Lee.
\newblock The moment formula for implied volatility at extreme strikes.
\newblock {\em Math. Finance}, 14(3):469--480, 2004.

\bibitem{Luc07}
V.~Lucic.
\newblock On singularities in the {H}eston model.
\newblock {\em Social Science Research Network. Preprint}, 2007.
\newblock http://papers.ssrn.com/sol3/papers.cfm?abstract\_id=1031222.

\bibitem{Luk70}
E.~Lukacs.
\newblock {\em Characteristic functions}.
\newblock Hafner Publishing Co., New York, 1970.
\newblock Second edition.

\bibitem{RevYor99}
D.~Revuz and M.~Yor.
\newblock {\em Continuous martingales and {B}rownian motion}, volume 293 of
  {\em Grundlehren der Mathematischen Wissenschaften [Fundamental Principles of
  Mathematical Sciences]}.
\newblock Springer-Verlag, Berlin, third edition, 1999.

\bibitem{Tit52}
E.~C. Titchmarsh.
\newblock {\em The Theory of Functions}.
\newblock Oxford University press, second edition, 1952.

\bibitem{Wid41}
D.~V. Widder.
\newblock {\em The {L}aplace {T}ransform}.
\newblock Princeton Mathematical Series, v. 6. Princeton University Press,
  Princeton, N. J., 1941.

\end{thebibliography}

\end{document}